\documentclass[12pt]{amsart}
\usepackage{a4wide}
\usepackage{amssymb}
\usepackage{amsthm}
\usepackage{amsmath}
\usepackage{amscd}
\usepackage{verbatim}

\numberwithin{equation}{section}

\theoremstyle{plain}
\newtheorem{theorem}{Theorem}[section]
\newtheorem{corollary}[theorem]{Corollary}
\newtheorem{lemma}[theorem]{Lemma}
\newtheorem{proposition}[theorem]{Proposition}

\newtheorem{conjecture}[theorem]{Conjecture}
\newtheorem{question}[theorem]{Question}

\theoremstyle{definition}
\newtheorem{definition}[theorem]{Definition}
\newtheorem{remark}[theorem]{Remark}

\theoremstyle{remark}

\newcommand{\CC}{\mathcal C \mathcal C}
\renewcommand{\O}{\mathcal O}

\newcommand{\R}{\mathbb{R}}
\newcommand{\Q}{\mathbb{Q}}
\newcommand{\Z}{\mathbb{Z}}
\newcommand{\N}{\mathbb{N}}
\newcommand{\C}{\mathbb{C}}

\renewcommand{\H}{\mathbb{H}}

\newcommand{\zxz}[4]{\begin{pmatrix} #1 & #2 \\ #3 & #4 \end{pmatrix}}

\newcommand{\kzxz}[4]{\left(\begin{smallmatrix} #1 & #2 \\ #3 & #4\end{smallmatrix}\right) }

\newcommand{\calC}{\mathcal{C}}

\newcommand{\calI}{\mathcal{I}}

\newcommand{\calK}{\mathcal{K}}

\newcommand{\calM}{\mathcal{M}}

\newcommand{\calO}{\mathcal{O}}

\newcommand{\frakd}{\mathfrak d}

\newcommand{\frakl}{\mathfrak l}

\newcommand{\CM}{\mathcal{CM}}
\newcommand{\eps}{\varepsilon}
\newcommand{\bs}{\backslash}
\newcommand{\norm}{\operatorname{N}}

\newcommand{\tr}{\operatorname{tr}}

\newcommand{\Sl}{\operatorname{SL}}

\newcommand{\Orth}{\operatorname{O}}

\newcommand{\Gal}{\operatorname{Gal}}

\newcommand{\CL}{\operatorname{CL}}

\newcommand{\ord}{\operatorname{ord}}

\newcommand{\dv}{\operatorname{div}}
\newcommand{\OK}{\mathcal{O}_K}

\newcommand{\darstell}[2]{#1 \mathcal{O}_F + #2 \mathcal{O}_F}

\begin{document}

\title[Twisted Borcherds products on Hilbert modular surfaces]{Twisted Borcherds products on Hilbert modular surfaces and their CM values}

\date{\today}
\author[Jan H.~Bruinier and Tonghai Yang]{Jan
Hendrik Bruinier and Tonghai Yang}
\address{Mathematisches Institut, Universit\"at zu K\"oln, Weyertal 86--90, D-50931 K\"oln, Germany}
\email{bruinier@math.uni-koeln.de }
\address{Department of Mathematics, University of Wisconsin Madison, Van Vleck Hall, Madison, WI 53706, USA}
\email{thyang@math.wisc.edu} \subjclass[2000]{11G15, 11F41, 14K22}
\thanks{ T.Y. Yang is partially supported by NSF grants DMS-0302043, 0354353, and
a NSA grant.}

\begin{abstract}
We construct a natural family of rational functions $\tilde\Psi_m$ on
a Hilbert modular surface from the classical $j$-invariant and its
Hecke translates.
These functions are obtained by means of a multiplicative
analogue of the Doi-Naganuma lifting and can be viewed as twisted
Borcherds products. We then study when the value of $\tilde\Psi_m$
at a CM point associated to a non-biquadratic quartic CM field
generates the `CM class field' of the reflex field. For the real
quadratic field $\Q(\sqrt{5})$,  we factorize the norm  of some of
these CM values  to $\Q(\sqrt 5)$ numerically.
\end{abstract}

\maketitle
\section{Introduction}
\label{sect:1}

In his papers \cite{Bo1} and \cite{Bo2} Borcherds constructed a
lifting from certain weakly holomorphic elliptic modular forms of
weight $1-n/2$ to meromorphic modular forms on the orthogonal group of
a rational quadratic space of signature $(2,n)$. Here we consider a
related construction in the particular case that $n=2$.

Let $p$ be a prime congruent to $1$ modulo $4$ and let
$F=\Q(\sqrt{p})$. We write $\calO_F$ for the ring of integers and
$\partial_F$ for the different of $F$.  Considering the lattice
$L^0=\Z^2\oplus \calO_F$ in the rational quadratic space
$L^0\otimes_\Q\Q$ of signature $(2,2)$, one obtains
from Borcherds' result a lifting from weakly holomorphic elliptic
modular forms of weight zero for the group $\Gamma_0(p)$ with
Nebentypus character $\epsilon_p=(\frac{\cdot}{p})$ to meromorphic
Hilbert modular forms for the Hilbert modular group
$\Gamma=\Sl_2(\calO_F)$ \cite{Br1,BB}.  This lifting can be viewed as
a multiplicative analogue of the Naganuma lift from holomorphic
modular forms of weight $k$ for $\Gamma_0(p)$ with Nebentypus
$\epsilon_p$ to holomorphic Hilbert modular forms of weight $k$ for
$\Gamma$ \cite{Na, Za}.

There is another lifting from holomorphic elliptic modular forms to Hilbert modular forms, namely the celebrated Doi-Naganuma lift which maps holomorphic modular forms of weight $k$ for $\Sl_2(\Z)$ to holomorphic Hilbert modular forms of weight $k$ for $\Gamma$ \cite{DN}. It was pointed out by Zagier that this lifting should have a multiplicative analogue as well \cite{Za3}. Moreover, Zagier stated several properties of such a multiplicative lifting and suggested that a proof could probably be given following the argument of \cite{Br1}.
One purpose of the present paper is to work out a proof along these lines.

Let $\H$ be the upper complex half plane, and put $q=e(\tau)=e^{2\pi
i\tau}$ for $\tau\in \H$.  Recall that a weakly holomorphic modular
form for $\Sl_2(\Z)$ is a meromorphic modular form for $\Sl_2(\Z)$
which is holomorphic outside the cusp $\infty$. In particular, every
weakly holomorphic modular form of weight zero for $\Sl_2(\Z)$ is a
polynomial in the $j$-function (which we normalize such that
$j(\tau)=q^{-1}+744+O(q)$).
%
In Section \ref{sect:5} we shall prove the following theorem. See
Theorem~\ref{hilbert} for a more detailed statement.

\begin{theorem}\label{hilbertintro}
Let $f=\sum_{n\gg -\infty}c(n)q^n\in \Z[j]$ be a weakly holomorphic
modular form of weight $0$ for $\Sl_2(\Z)$ with integral Fourier
coefficients. Then there exists a symmetric meromorphic Hilbert
modular function $\Psi(z,f)$ for $\Gamma$ (of weight $0$, with trivial
multiplier system, defined over $F$) such that:
\begin{enumerate}
\item[(i)] The divisor of $\Psi(z,f)$ is determined by the polar
part of $f$ at the cusp $\infty$. It equals $\sum_{n>0}  c(-n) \tilde T_{n}$,
where $\tilde T_{n}$ denotes the twisted Hirzebruch-Zagier divisor of discriminant $n$ defined in Section \ref{sect:3}.
\item[(ii)]
The function $\Psi(z,f)$ has the Borcherds product expansion
\begin{align*}
\Psi(z,f)&= \prod_{\substack{\nu\in\partial_F^{-1} \\
\nu>0}}\prod_{b\,(p)} \left(1-e(\tfrac{b}{p}+\nu
z_1+\nu'z_2)\right)^{\epsilon_p(b)c(p\nu\nu')},
\end{align*}
which converges normally for all $z=(z_1,z_2)\in \H^2$ with
$\Im(z_1)\Im(z_2)> Np$ outside the set of poles, where $N=\max
\{n\in \Z;\; c(-n)\neq 0\}$. In particular, $\Psi(z, C) =1$ for
any constant $C$.
\item[(iii)] The lifting is multiplicative, i.e., if $f,g\in \Z[j]$,
then $\Psi(f+g)=\Psi(f)\Psi(g)$.
\end{enumerate}
\end{theorem}

We now briefly describe the idea of the proof.
For a positive integer $m$, we define in Section
\ref{sect:3} a certain ``twisted Hirzebruch-Zagier divisor'' $\tilde
T_m$ on the Hilbert modular surface corresponding to $\Gamma$ (see
also \cite{Za2} pp.~166).  In Section \ref{sect:4}, following
\cite{Br1}, we construct an automorphic Green function $\tilde
\Phi_{mp^2}(z,s)$ for  $\tilde T_m$.  We study its main
properties and compute its Fourier expansion (Theorems
\ref{fourier1} and \ref{fourier2}).
By means of an identity relating certain finite exponential sums to
Kloosterman sums (Lemma \ref{GaHb}), we find that the Fourier
coefficients of $\tilde \Phi_{mp^2}(z,s)$ are closely related to the
coefficients of non-holomorphic Poincar\'e series of weight zero for
$\Sl_2(\Z)$ (see Section \ref{sect:2}).  Using the fact that any
weakly holomorphic modular form of weight zero can be uniquely written
as a linear combination of the non-holomorphic Poincar\'e series,
Theorem \ref{hilbertintro} can be deduced.  Here the main point is an
identity expressing $\log|\Psi(z,f)|$ as a linear combination of the
automorphic Green functions $\tilde\Phi_{mp^2}(z,1)$, see Theorem
\ref{hilbert} (iv).

An alternative proof could be given by interpreting the lifting as a
regularized theta lifting for the dual pair $\Sl_2(\R)$, $\Orth(2,2)$
as in \cite{Bo2,Br2}, and by considering suitable ``twists'' of Siegel
theta functions as kernel functions.
It would actually be very interesting to describe
such twisted Borcherds liftings in greater generality for
$\Orth(2,n)$.  However, it seems not quite clear what the right
``twists'' of Siegel theta functions should be in general.  We have
not pursued this approach in the present paper, because the proof
using automorphic Green functions leads to the result in a direct way.
Moreover, the construction of the Green functions
should be of independent interest.

If $m$ is a positive integer, we write $\tilde\Psi_m$ for the symmetric
Hilbert modular function of weight $0$ which is the twisted Borcherds
lift via Theorem \ref{hilbertintro} of the unique weakly
holomorphic modular form
\[
J_m=q^{-m}+O(q)\in \Z[j].
\]
In this way we obtain a `canonical' family of rational functions
on the Hilbert modular surface associated to $\Gamma$. It seems
natural to ask for their arithmetic properties. For instance, one
can ask whether $\tilde \Psi_1$, the lifting of the $j$-function,
has a nice moduli interpretation as an invariant of abelian
surfaces with $\calO_F$-multiplication and
$\partial_F^{-1}$-polarization.

It is a well-known and beautiful fact that $j(\frac{D+\sqrt{-D}}2)$
generates the Hilbert class field of the imaginary quadratic field
$\Q(\sqrt{-D})$. 
In Section \ref{sect:class}, we study when the value of
$\tilde\Psi_m$ at a CM point $z$ associated to a quartic CM number
field $K$ with totally real subfield $F$
generates the `CM class field' $H_{\tilde K}$ of the reflex field
$\tilde K$ (see Section 6 for precise definitions).
 In particular,
we show that it generates this CM class field when the CM value $\tilde\Psi_m(\CM(K))$, the evaluation of
$\tilde\Psi_m$ at the CM cycle $\CM(K)$ corresponding to $K$, is not an odd power in $F$ (Corollary
\ref{clcor1.3}). We also prove the following theorem (see Theorem \ref{cltheo6.6} for details).

\begin{theorem} Let $F=\mathbb Q(\sqrt p)$ be a fixed real
quadratic field with $p \equiv 1 \mod 4$ prime. Then there is a
constant $d >0$ such that for any CM quartic field $K$ of
discriminant $d_K=p^2 q$ with $q \equiv 1 \mod 4$ prime, and a CM
point $z$ in the Hilbert modular surface $X=\hbox{SL}_2(\O_F)
\backslash \mathbb H^2$ of CM type $(K, \Phi)$ by $\O_K$, one has:
\begin{enumerate}
\item[(i)] $M(\tilde\Psi_1(z), \tilde\Psi_2(z))$ is an unramified abelian extension of $M$, where $M$ is the
smallest Galois extension of $\mathbb Q$ containing $K$.
\item[(ii)]
Let $H_{\tilde K}$ be the `CM class field' of the reflex field $\tilde K$, and let $L_K=M H_{\tilde K}$. Then
$M(\tilde\Psi_1(z), \tilde\Psi_2(z))$ is a subfield of $L_K$ with bounded index $[L_K:M(\tilde\Psi_1(z),
\tilde\Psi_2(z))] \le d$.
\item[(iii)]
One has $\lim_{q \rightarrow \infty} \frac{\log
[M(\tilde\Psi_1(z), \tilde\Psi_2(z)): M]}{\log  \sqrt q}=1$.
\end{enumerate}
\end{theorem}

In \cite{BY} the authors derived a formula for the values of
(untwisted) Borcherds products (in the sense of \cite{BB}) at CM
cycles $\CM(K)$.
It would be very interesting to obtain an analogous formula for the CM
values of the twisted Borcherds products of Theorem
\ref{hilbertintro}. In fact one can ask if it is possible to modify
the proof of \cite{BY} to give such a result. One key ingredient of
the proof -- the relation between Borcherds products and automorphic
Green functions -- is already worked out in the present
paper. However, in the second main step of the argument
it is not clear at all how the function $\psi$ on the lattice $L^0$
defined in \eqref{eq:psi} translates to some natural function on the
reflex field $\tilde K$.

In Section \ref{sect:examples} we study some examples in the special case that $F=\Q(\sqrt{5})$. We write $\tilde \Psi_1$ and $\tilde \Psi_2$ in terms of the generators of the ring of symmetric Hilbert modular forms of even weight given by Gundlach \cite{Gu}. This can be used to compute the Fourier expansions explicitly, which in turn can be employed to compute some CM values (mainly) numerically.
For instance, for the CM point $z_0=(\zeta_5,\zeta_5^2)$ (where $\zeta_5=e^{2\pi i/5}$) corresponding to the
cyclic CM extension $K=\Q(\zeta_5)$ of $F=\Q(\sqrt{5})$ we find that
\begin{align*}
\tilde \Psi_1(z_0) &=\frac{156973921227+70200871784\sqrt{5}}{156973921227-70200871784\sqrt{5}}=
\frac{\omega^{27}\cdot (4+\omega')^5\cdot(5+\omega')^5}{\omega'{}^{27}\cdot(4+\omega)^5\cdot(5+\omega)^5},
\end{align*}
where $\omega=\frac{1+\sqrt{5}}{2}$.
By means of the results of \cite{BY} we derive a heuristic how the CM values of twisted Borcherds products should look like. In particular, we obtain a conjecture on the prime ideals  $\frakl\subset\calO_F$ at which the CM value of a twisted Borcherds product can have non-zero order (see Conjecture \ref{conj}). The same phenomenon as in \cite{GZ} and \cite{BY} should happen: Such prime ideals should be of small norm.

Finally, in Section \ref{sect:7}, we list some open problems for
further research. For instance, in all the examples we computed it
turned out that $\tilde \Psi_1(\CM(K, \Phi, \O_F))$ belongs to the
field $F$. Moreover, if $K/\Q$ is non-Galois then $\tilde
\Psi_1(\CM(K, \Phi, \O_F))$ is square-free. According to
Corollary~\ref{clcor1.3}, this implies that for $z\in \CM(K, \Phi,
\O_F)$ the CM value $\tilde \Psi_1(z)$ generates the class field
$L_K$ over $M$. It is an interesting question whether this is a
general phenomenon.

  We mention that J. Rouse has used Theorem 1.1 to determine the Fourier
coefficients of modular functions $f \in \mathbb{Z}[j]$ in terms
of traces of singular moduli
  \cite{Rouse}.


\medskip

We thank Heike Hippauf and Sebastian Mayer for their valuable help
with the computations in Section \ref{sect:examples}. Moreover, we
thank  Lev Borisov, Bas Edixhoven,  and  Don Zagier  for
interesting and useful discussions, and thank Shou-Wu Zhang for
bringing to our attention his equidistribution theorem, which is
needed in the proof of Theorem \ref{cltheo6.6}.

\section{Non-holomorphic Poincar\'e series}
\label{sect:2}

Here we consider non-holomorphic Poincare series of weight $0$. The
results of this section are known. We state them for completeness and
to fix the notation. For details we refer to \cite{He}, \cite{Ni},
\cite{Br2}.

Let $I_\nu(z)$ and $K_\nu(z)$ be the usual modified Bessel functions as in \cite{AS} \S10. For convenience we put
for $s\in\C$ and $y\in\R\setminus \{0\}$:
\begin{align}\label{calI}
\calI_s(y)&=  
\sqrt{\frac{\pi |y|}{2}}I_{s-1/2}(|y|),\\
\calK_s(y)&= \label{calK}
\sqrt{\frac{2|y|}{\pi}}K_{s-1/2}(|y|).
\end{align}
The functions $\calI_s(y)$ and  $\calK_s(y)$ are holomorphic in $s$. At $s=1$ they have the special values
\begin{align}
\label{Mspecial}
\calI_{1}(y)&=\sinh(|y|), \\
\label{Wspecial}
\calK_{1}(y)&=e^{-|y|},\\
\label{MWspecial} 2\calI_{1}(y)+\calK_{1}(y)&=e^{|y|}.
\end{align}

The full elliptic modular group $\Gamma'=\Sl_2(\Z)$ acts on the upper complex half plane $\H=\{\tau\in \C;\;\Im(\tau)>0\}$ by linear fractional transformations. We write $\Gamma_\infty'=\{\kzxz{1}{n}{0}{1};\;n\in \Z\}$. As usual we abbreviate $e(x)=e^{2\pi i x}$.

For a positive integer $m$
we define the Poincar\'e series of weight $0$ and index $m$ by
\begin{equation}\label{DefF}
F_{m}(\tau,s)=\sum_{\gamma\in \Gamma'_\infty\bs \Gamma'} \calI_s\big(2\pi m \Im(\gamma \tau)\big) e\big(-m \Re(\gamma \tau)\big),
\end{equation}
where $\tau=x+iy\in\H$ and $s\in\C$ with $\Re(s)>1$. It converges normally for $\Re(s)>1$ and defines a
$\Gamma'$-invariant function on $\H$. It is an eigenfunction of the hyperbolic Laplacian with eigenvalue $s(s-1)$.

\begin{theorem}\label{fourierF}
The Poincar\'e series $F_{m}(\tau,s)$ has the Fourier expansion
\begin{align*}
F_{m}(\tau,s)&=\left(2\calI_{s}(2\pi m y)+\calK_s(2\pi m y)\right) e(-mx)\\
&\phantom{=}{}+ b_m(0,s) y^{1-s}  + \sum_{\substack{ n\in \Z\setminus \{0\}}} b_m(n,s) \calK_s(2\pi ny) e(nx),
\end{align*}
where
\[
b_m(n,s)=\begin{cases}
\displaystyle 2\pi \left| \frac{m}{n}\right|^{1/2} \sum_{c=1}^\infty  H_c(m,n) I_{2s-1}\!\left(\frac{4\pi}{c}\sqrt{|mn|}\right), & \text{$n>0$,}\\[3ex]
\displaystyle \frac{4 \pi^{1+s} m^{s}}{(2s-1)\Gamma(s)} \sum_{c=1}^\infty c^{1-2s} H_c(m,0), & \text{$n=0$,}\\[3ex]
\displaystyle -\delta_{-m,n} + 2\pi \left| \frac{m}{n}\right|^{1/2} \sum_{c=1}^\infty H_c(m,n)
J_{2s-1}\!\left(\frac{4\pi}{c}\sqrt{|mn|}\right), & \text{$n<0$}.
\end{cases}
\]
Here $ H_c(m,n)$ denotes the Kloosterman sum
\begin{equation}\label{DefH_c}
H_c(m,n)=\frac{1}{c} \sum_{\substack{d(c)^*}} e\!\left( \frac{nd-m\bar d}{c}\right),
\end{equation}
where the sum runs through the multiplicative group $(\Z/c\Z)^*$
and $\bar d$ denotes the multiplicative inverse of $d$. Moreover,
$J_\nu(z)$ and $I_\nu(z)$ are the usual Bessel functions as
defined in \cite{AS}~\S9.
\end{theorem}

\begin{proof}
This is a special case of \cite{Br2} Theorem 1.9. See also \cite{Ni} or \cite{He}.
\end{proof}

Notice that $H_c(m,n)=H_c(n,m)$.

\begin{proposition}\label{constantterm}
The constant term of $F_m(\tau,s)$ is equal to
\[
b_m(0,s)=\frac{4 \pi}{(2s-1)} \frac{\sigma_m(2s-1)}{\pi^{-s}\Gamma(s)\zeta(2s)}.
\]
Here $\zeta(s)$ denotes the Riemann zeta function and $\sigma_m(s)$ the divisor sum
\[
\sigma_m(s)=m^{(1-s)/2}\sum_{d\mid m} d^s.
\]
\end{proposition}

\begin{proof}
By Theorem \ref{fourierF} we have
\begin{align*}
b_m(0,s)&= \frac{4 \pi^{1+s} m^{s}}{(2s-1)\Gamma(s)} \sum_{c=1}^\infty c^{-2s} \sum_{d\,(c)^*}
e\left(\frac{md}{c}\right).
\end{align*}
If we insert the formula for the Ramanujan sum (see \cite{Ap} Chapter 8.3),
\[
\sum_{d\,(c)^*} e\left(\frac{md}{c}\right) = \sum_{a\mid (c,m)}\mu(c/a)a,
\]
where $\mu$ is the Moebius function, we obtain
\begin{align*}
b_m(0,s)&=\frac{4 \pi^{1+s} m^{s}}{(2s-1)\Gamma(s)} \sum_{a\mid m} \sum_{\substack{c=1\\a\mid c}}^\infty c^{-2s}\mu(c/a)a\\
&=\frac{4 \pi^{1+s}}{(2s-1)\Gamma(s)\zeta(2s)}  m^{s} \sum_{a\mid m} a^{1-2s}\\
&=\frac{4 \pi^{1+s}\sigma_m(2s-1)}{(2s-1)\Gamma(s)\zeta(2s)}.
\end{align*}
This proves the Proposition.
\end{proof}

Recall that a weakly holomorphic modular form for $\Gamma'$ is a
meromorphic modular form for $\Gamma'$ which is holomorphic outside
the cusp $\infty$. In particular, the space of weakly holomorphic
modular forms for $\Gamma'$ of weight $0$ is $\C[j]$, the polynomial
ring in the $j$-function.

\begin{theorem}\label{fourierF1}
The special value $F_{m}(\tau,1)$ has the Fourier expansion
\begin{align*}
F_{m}(\tau,1)&=q^{-m}
 +\sum_{n\geq 0} b_m(n,1) q^n,
\end{align*}
where $q=e^{2\pi i\tau}$ as usual. In particular, $F_{m}(\tau,1)$ is the unique weakly holomorphic modular form of
weight $0$ for $\Gamma'$ whose Fourier expansion starts $q^{-m}+b_m(0,1)+O(q)$.
\end{theorem}

\begin{proof}
If we insert the special values \eqref{Mspecial} and \eqref{Wspecial} in the Fourier expansion given in Theorem
\ref{fourierF1}, we find
\[
F_{m}(\tau,1)= q^{-m}+\sum_{n\geq 0} b_m(n,1) q^n + \sum_{\substack{ n<0}} b_m(n,1) e(n\bar\tau).
\]
This implies that
\[
\overline{\tfrac{\partial}{\partial\bar\tau} F_{m}(\tau,1)}=-2\pi i \sum_{\substack{ n<0}} \overline{b_m(n,1)} n
e(-n\tau)
\]
is a holomorphic modular form of weight $2$ for $\Gamma'$ and therefore has to vanish identically. We obtain the
assertion.
\end{proof}

\begin{remark}
We have $b_m(0,1)=24\sigma_m(1)$.
\end{remark}

\section{Hilbert modular surfaces and Hirzebruch-Zagier divisors}
\label{sect:3}

Let $p\equiv 1 \pmod{4}$ be a prime and consider the real quadratic
field $F=\Q(\sqrt{p})$. We write $\calO_F$ for the ring of integers in
$F$ and $\partial_F=(\sqrt{p})$ for the different ideal. Moreover, we
write $\epsilon_p(x)=(\frac{x}{p})$ for the quadratic character of
$F$. The conjugation in $F$ is denoted by $x\mapsto x'$ and the norm
of $x\in F$ by $\norm(x)=xx'$.

The Hilbert modular group $\Gamma=\Sl_2(\calO_F)$ acts on $\H^2$ in the usual way. We denote by $X=\Gamma\bs \H^2$
the Hilbert modular surface associated to $F$.    We will use $z=(z_1,z_2)$ as a standard variable on $\H^2$ and
write $z_1=x_1+iy_1$, $z_2=x_2+iy_2$ for the decomposition in real and imaginary parts. Recall that a Hilbert
modular form $H(z_1,z_2)$ is called symmetric if $H(z_1,z_2)=H(z_2,z_1)$.  It is called skew-symmetric if
$H(z_1,z_2)=-H(z_2,z_1)$.

It is well known that the Hilbert
modular group can also be viewed as a discrete subgroup of the orthogonal group of the rational quadratic space
\begin{equation}
V=\{M\in \operatorname{Mat}_2(F);\; {}^tM=M'\}
= \{\kzxz{a}{\lambda}{\lambda'}{b};\; a,b\in \Q,\; \lambda\in F\},
\end{equation}
equipped with the quadratic form $M\mapsto \det(M)$. Here ${}^tM$ is the transpose of $M$. The group
$\Sl_2(F)$ acts on $V$ via
\begin{equation}
\gamma . M = \gamma' M \, {}^t\gamma
\end{equation}
leaving the quadratic form invariant. The lattice
\begin{align}\label{def:L}
L=\{\kzxz{a}{\lambda}{\lambda'}{b}\in V;\; a,b\in \Z,\; \lambda\in
\partial^{-1}_F\}
\end{align}
is stable under the action of $\Gamma$. Recall that $M\in L$ is called {\em primitive}, if $\frac{1}{d} M\notin L$ for every integer $d>1$.

If $M=\kzxz{a}{\lambda}{\lambda'}{b}$ is a vector in $L$, we
let
\[
M^\perp=\{(z_1,z_2)\in \H^2; \quad az_1 z_2 +\lambda z_1 + \lambda' z_2 +b=0\}
\]
be the corresponding analytic divisor on $\H^2$.
If $m$ is a positive integer, then
\begin{align}
F_m&=\sum_{\substack{\text{$M\in L/\{\pm 1\}$ primitive}\\
\det(M)=m/p}} M^\perp
\end{align}
is a $\Gamma$-invariant divisor on $\H^2$, which descends to an algebraic divisor on $X$, also denoted by $F_m$.
It is well known that the Hirzebruch-Zagier divisors $T_m$ on $X$ can be written as $T_m=\sum_{d^2|m} F_{m/d^2}$ (see \cite{Ge}, \cite{HZ}).

The divisor $F_m$  on $X$ is irreducible if and only if $p^2\nmid m$. If $p^2|m$ then it decomposes into two irreducible
components $F_m=F_m^++F_m^-$. To distinguish these components we define a function on lattice vectors
$M=\kzxz{a}{\lambda}{\lambda'}{b}\in L$ of norm $ab-\lambda\lambda'\in p\Z$ divisible by $p$ by
\begin{align}\label{eq:psi}
\psi(M)=\begin{cases}
 (\tfrac{a}{p}),& \text{if $p\nmid a$,}\\
 (\tfrac{b}{p}),& \text{if $p\nmid b$,}\\
0,&\text{if $p\mid (a,b)$.}
            \end{cases}
\end{align}
The following lemma shows that $\psi$ is well defined.

\begin{lemma}
Let $M=\kzxz{a}{\lambda}{\lambda'}{b}\in L$ be a $p$-primitive vector (i.e. $\frac{1}{p}M\notin L$) and assume that
$ab-\lambda\lambda'=np$ for some $n\in \Z$. Then $(\tfrac{a}{p})+(\tfrac{b}{p})\neq 0$.
\end{lemma}

\begin{proof}
First we notice that the hypothesis implies $p\nmid (a,b)$. This immediately yields the assertion if $p\mid a$, or
$p\mid b$, respectively. So we may assume that $p\nmid a$ and $p\nmid b$. We have to show that
$(\tfrac{a}{p})=(\tfrac{b}{p})$. If we write $\lambda=\frac{1}{2}(c+d\sqrt{p})$ with $c,d\in
\Z$ we find that
\begin{align*}
4ab=c^2-pd^2+4np\equiv \alpha^2\pmod{p}.
\end{align*}
This implies the assertion.
\end{proof}

The function $\psi$ has the following genus character interpretation.
Assume that  $M=\kzxz{a}{\lambda}{\lambda'}{b}\in L$ and
$ab-\lambda\lambda'\in p\Z$.
Then $\lambda$ actually belongs to $\calO_F$ and we may write $\lambda=\frac{1}{2}(c+d\sqrt{p})$ with $c,d\in
\Z$. We have
\begin{align}\label{quadform}
\zxz{a}{\lambda}{\lambda'}{b}\equiv \zxz{a}{c/2}{c/2}{b}\pmod{\partial_F}
\end{align}
and the latter matrix defines a binary integral quadratic form $Q=[a,c,b]$ of discriminant $c^2-4ab$ divisible by $p$. Recall that on  binary integral quadratic forms of discriminant divisible by $p$ we have a genus character $\chi_p$ which is defined as follows:
\[
\chi_p(Q)=\begin{cases} (\frac{n}{p}),& \text{if $p\nmid Q$ and $(n,p)=1$ and $Q$ represents $n$},\\
0,& \text{if $p\mid Q$.}
      \end{cases}
\]
This definition does not depend on the choice of $n$. If the quadratic form $Q$ corresponds to $M$ as in \eqref{quadform},
then $\psi(M)=\chi_p(Q)$.

One easily verifies that the value of $\psi(M)$ only depends on the $\Gamma$-orbit of $M$.
In fact, it is well known that the components of $F_{mp^2}$ are distinguished by the values of the function $\psi$ (see
\cite{Za2}, \cite{Ge} Chapter V.3).


\begin{definition}
For a positive integer $m$ we define the twisted Hirzebruch-Zagier divisor of index $m$ by
\begin{align}\label{def: tildetm}
\tilde T_m = \sum_{\substack{M\in L/\{\pm 1\}\\
\det(M)=mp}} \psi(M)\cdot M^\perp.
\end{align}
\end{definition}

For instance, if $m$ is square-free then $\tilde T_m  = F_{mp^2}^+-F_{mp^2}^-$.

\section{Automorphic Green functions}
\label{sect:4}

Let $m$ be a positive integer.
Following \cite{Br1} we define the automorphic Green function corresponding to $\tilde T_m$ by
\begin{align}
\tilde \Phi_{mp^2}(z_1,z_2,s)=
\sum_{\substack{\kzxz{a}{\lambda}{\lambda'}{b} \in L\\
ab-\norm(\lambda)=mp}} \psi\kzxz{a}{\lambda}{\lambda'}{b} Q_{s-1}\left(1+\frac{ |az_1
z_2 +\lambda z_1+\lambda' z_2+b|^2}{ 2y_1 y_2 mp }\right),
\end{align}
where $Q_{s-1}(t)$ is the Legendre function of the second kind defined by (cf.~\cite{AS} \S8)
\[
Q_{s-1}(t) =\int_0^\infty (t +\sqrt{t^2-1} \cosh v)^{-s} dv\qquad (t>1, \,\Re(s)>0).
\]

 The sum
converges normally for $(z_1,z_2)\in  \H^2\setminus \tilde T_m$ and $s\in \C$ with $\Re(s)>1$. We will continue it
to a neighborhood of $s=1$ by computing the Fourier expansion. To this end we write
\begin{align}\label{auft}
\tilde\Phi_{mp^2}(z_1,z_2,s)&= \tilde\Phi^0_{mp^2}(z_1,z_2,s) +
2\sum_{a=1}^{\infty} \tilde\Phi^a_{mp^2}(z_1,z_2,s)\\ \intertext{with}
\tilde\Phi^a_{mp^2}(z_1,z_2,s)&= \sum_{\substack{ b\in\Z \\
\lambda\in\partial_F^{-1}\\ ab-\norm(\lambda)=mp }}
\psi\kzxz{a}{\lambda}{\lambda'}{b} Q_{s-1}\left(1+\frac{ |az_1 z_2
+\lambda z_1+\lambda' z_2+b|^2}{ 2y_1 y_2 mp }\right). \label{Phia}
\end{align}
Note that the partial sums $\tilde\Phi^a_{mp^2}(z_1,z_2,s)$ converge normally for $\Re(s)>1/2$.

For $a=0$ we have
\[
\tilde\Phi^0_{mp^2}(z_1,z_2,s)= \sum_{\substack{ b\in\Z \\ \lambda\in\partial_F^{-1}\\ \norm(\lambda)=-mp }}
\epsilon_p(b)  Q_{s-1}\left(1+\frac{ |\lambda z_1+\lambda' z_2+b|^2}{ 2y_1 y_2 mp }\right).
\]
If $a$ is a positive integer coprime to $p$, then
\begin{align*}
\tilde\Phi^a_{mp^2}(z_1,z_2,s)&= \epsilon_p(a)  \sum_{\substack{ b\in\Z \\ \lambda\in\partial_F^{-1}\\ ab-\norm(\lambda)=mp }} Q_{s-1}\left(1+\frac{ |az_1 z_2 +\lambda z_1+\lambda' z_2+b|^2}{ 2y_1 y_2 mp }\right)\\
&= \epsilon_p(a) \Phi^a_{mp^2}(z_1,z_2,s),
\end{align*}
where $\Phi^a_{mp^2}(z_1,z_2,s)$ is the function defined in \cite{Br1} (see \S3 equation (16)). If $a$ is
divisible by $p$, then
\begin{align*}
\tilde\Phi^a_{mp^2}(z_1,z_2,s)&= \sum_{\substack{ b\in\Z \\ \lambda\in\partial_F^{-1}\\ ab-\norm(\lambda)=mp }}
\epsilon_p(b)  Q_{s-1}\left(1+\frac{ |az_1 z_2 +\lambda z_1+\lambda' z_2+b|^2}{ 2y_1 y_2 mp }\right).
\end{align*}

We compute the Fourier expansion of $\tilde\Phi^a_{mp^2}(z_1,z_2,s)$ in these  three cases. We put for a positive
real number $A$,
\begin{equation}\label{HAs}
H^A_s(z_1,z_2) = \sum_{\theta\in\calO}  Q_{s-1}\left(1+\frac{ |(z_1+\theta)(z_2+\theta')+A|^2  }{2y_1 y_2
A}\right),
\end{equation}
and denote (for $y_1y_2>A$) the Fourier expansion by
\[
H^A_s(z_1,z_2) = \sum_{\nu\in\frakd^{-1}} b_s^A(\nu,y_1,y_2) e(\nu x_1+\nu'x_2).
\]
Moreover, we let $R(mp^2)$ be a set of representatives for
\begin{align*}
&\{\lambda\in \partial_F^{-1}/a\calO_F;\quad \norm(\lambda\sqrt{p})\equiv mp^2 \pmod{ap}\}\\
&=\{\lambda\in \calO_F/a\calO_F;\quad \norm(\lambda)\equiv -mp \pmod{a}\}.
\end{align*}

We start with the case that $a$ is positive and coprime to $p$. Here we can argue as in  \cite{Br1}. We may write
\begin{align*}
\tilde\Phi^a_{mp^2} (z_1,z_2,s)&= \epsilon_p(a)\sum_{\lambda\in R(mp^2)}\sum_{\theta\in\calO} Q_{s-1}\left(1+\frac{ |(z_1+\theta+\tfrac{\lambda'}{a})(z_2+\theta'+\tfrac{\lambda}{a})+\tfrac{mp}{a^2}|^2 }{ 2y_1 y_2 mp/a^2}\right)\\
&=\sum_{\nu\in\partial_F^{-1}}\epsilon_p(a)\sum_{\lambda\in R(mp^2)} e\left(\frac{\tr(\nu\lambda)}{
a}\right)b_s^{mp/a^2}(\nu,y_1,y_2) e(\nu z_1+\nu'z_2).
\end{align*}
When $a$ is positive and divisible by $p$, one finds in a similar way that
\begin{align*}
\tilde\Phi^a_{mp^2} (z_1,z_2,s)&= \sum_{\lambda\in R(mp^2)}\sum_{\theta\in\calO} \epsilon_p(\tfrac{\lambda\lambda'+mp}{a}) Q_{s-1}\left(1+\frac{ |(z_1+\theta+\tfrac{\lambda'}{a})(z_2+\theta'+\tfrac{\lambda}{a})+\tfrac{mp}{a^2}|^2 }{ 2y_1 y_2 mp/a^2}\right)\\
&=\sum_{\nu\in\partial_F^{-1}}
\sum_{\lambda\in R(mp^2)} \epsilon_p(\tfrac{\lambda\lambda'+mp}{a})  e\left(\frac{\tr(\nu\lambda)}{
a}\right)b_s^{mp/a^2}(\nu,y_1,y_2) e(\nu z_1+\nu'z_2).
\end{align*}
If we define
\begin{equation}\label{tildeGa}
\tilde G_a(mp^2,\nu)=\begin{cases} \epsilon_p(a)\sum_{\lambda\in R(mp^2)}
e\left(\frac{\tr(\nu\lambda)}{ a}\right),&\text{if $p\nmid a$,}\\[1ex]
\sum_{\lambda\in R(mp^2)} \epsilon_p(\tfrac{\lambda\lambda'+mp}{a}) e\left(\frac{\tr(\nu\lambda)}{
a}\right),&\text{if $p\mid a$,}
             \end{cases}
\end{equation}
we may finally write
\begin{equation}\label{Phim2}
\tilde\Phi_{mp^2} (z_1,z_2,s) = \tilde\Phi^0_{mp^2}(z_1,z_2,s) + 2
\sum_{\nu\in\partial_F^{-1}}\left[\sum_{a=1}^{\infty} \tilde G_a(mp^2,\nu)b_s^{mp/a^2}(\nu,y_1,y_2) \right]e(\nu
x_1+\nu'x_2).
\end{equation}

For $r_1,r_2\in \R$ we briefly write
\begin{align*}
\alpha(r_1,r_2)&:= \max(|r_1|,|r_2|),\\
\beta(r_1,r_2)&:= \min(|r_1|,|r_2|).
\end{align*}

\begin{lemma}\label{lphi0}
The function $\tilde\Phi^0_{mp^2}(z_1,z_2,s)$ has the Fourier expansion
\begin{multline}\label{lphi0a}
\tilde\Phi_{mp^2}^0(z_1,z_2,s) =  2 \sqrt{p} \sum_{\substack{ \lambda\in\partial_F^{-1}\\ \norm(\lambda)=-m/p}}
\sum_{n\geq 1} \frac{\epsilon_p(n)}{n}
 \calI_{s}(2\pi n \beta(\lambda y_1,\lambda' y_2))\\
{}\times \calK_{s}(2\pi n \alpha(\lambda y_1,\lambda' y_2))e(n\lambda x_1+n\lambda' x_2).
\end{multline}
\end{lemma}

\begin{proof}
By definition we have
\begin{align*}
\tilde\Phi_{mp^2}^0(z_1,z_2,s) &= \sum_{\substack{b\in \Z,\,\lambda\in \partial_F\\ \lambda\lambda'=-mp}}
\epsilon_p(b)
Q_{s-1}\left(1+\frac{ |\lambda z_1+\lambda' z_2+b|^2}{ 2y_1 y_2 mp }\right)\\
&= \sum_{\substack{b\,(p),\,\lambda\in \partial_F^{-1}\\ \lambda\lambda'=-m/p}} \epsilon_p(b) \sum_{b'\in \Z}
Q_{s-1}\left(1+\frac{ |\lambda z_1+\lambda' z_2+b'+b/p|^2}{ 2y_1 y_2 m/p }\right).
\end{align*}
We find that
\begin{align}\label{Phi0h}
\tilde\Phi_{mp^2}^0(z_1,z_2,s) = \sum_{\substack{b\,(p)}} \sum_{\substack{\lambda\in \partial_F^{-1}\\
\lambda\lambda'=-m/p}}\epsilon_p(b) h_{\alpha(\lambda y_1, \lambda' y_2), \beta(\lambda y_1, \lambda'
y_2)}(\lambda x_1+ \lambda' x_2+b/p),
\end{align}
where
\[
h_{\alpha,\beta}(x) = \sum_{b\in\Z}  Q_{s-1}\left(\frac{(x+b)^2+\alpha^2+\beta^2}{2\alpha\beta}\right).
\]
By \cite{Br1} Lemma 1, for $\alpha>\beta>0$, the function $h_{\alpha,\beta}(x)$ has the Fourier expansion
\[
h_{\alpha,\beta}(x)= \frac{2\pi}{2s-1} \alpha^{1-s}\beta^{s} +\sum_{n\in\Z-\{0\}}\frac{1}{|n|} \calI_{s}(2\pi n
\beta) \calK_{s}(2\pi n\alpha) e(nx).
\]
Inserting this into \eqref{Phi0h}, we obtain the assertion.
\end{proof}

We recall the following lemma from \cite{Br1}:

\begin{lemma}\label{lHA}
Let $y_1 y_2>A>0$. The function $H^A_s(z_1,z_2)$ defined by (\ref{HAs}) has the Fourier expansion
\begin{equation*}
H^A_s(z_1,z_2) = \sum_{\nu\in\partial_F^{-1}} b_s^A(\nu,y_1,y_2) e(\nu x_1+\nu'x_2)
\end{equation*}
with
\begin{align*}
b^A_s(0,y_1,y_2) &= \frac{\pi \Gamma(s-1/2)^2}{ 2\sqrt{p}\,\Gamma(2s)}(4A)^s (y_1 y_2)^{1-s},\\
b^A_s(\nu,y_1,y_2)&=\pi\sqrt{\frac{A}{ p|\nu\nu'|}}\,I_{2s-1}(4\pi\sqrt{A|\nu\nu'|})\calK_{s}(2\pi |\nu| y_1) \calK_{s}(2\pi |\nu'| y_2), \quad\nu\nu'>0,\\
b^A_s(\nu,y_1,y_2)&=\pi\sqrt{\frac{A}{ p|\nu\nu'|}}\,J_{2s-1}(4\pi\sqrt{A|\nu\nu'|})\calK_{s}(2\pi |\nu| y_1)
\calK_{s}(2\pi |\nu'| y_2), \quad\nu\nu'<0.
\end{align*}

\end{lemma}

The following identity of finite exponential sums is crucial for the main result of the present paper.

\begin{lemma}\label{GaHb}
Let $a\in\N$, $m\in\Z$, and $\nu\in\partial_F^{-1}$. Then
\begin{equation}\label{eq:key}
\frac{1}{ a } \,\tilde G_a(mp^2,\nu) = \sum_{\substack{r|\nu \\ r|a}} \epsilon_p(r) H_{a/r}\left( \frac{p\nu\nu'}{
r^2},m\right),
\end{equation}
where the finite exponential sums $\tilde G_a(m,\nu)$ resp.~$H_b(m,n)$ are defined by (\ref{tildeGa})
resp.~(\ref{DefH_c}).
\end{lemma}

\begin{proof}
We follow the proof of the Proposition in \cite{Za} \S4. Both sides in \eqref{eq:key} are clearly periodic in $m$
with period $a$. Therefore it suffices to show that the finite Fourier transforms are equal, i.e., that for every
$h\pmod{a}$ we have
\begin{align}\label{eq:key2}
\frac{1}{ a }\sum_{m\,(a)}e\left(-\frac{hm}{a}\right) \tilde G_a(mp^2,\nu) =
\sum_{m\,(a)}e\left(-\frac{hm}{a}\right)\sum_{\substack{r|\nu \\ r|a}} \epsilon_p(r) H_{a/r}\left(
\frac{p\nu\nu'}{ r^2},m\right).
\end{align}
Inserting the definition of $H_{c}\left( n,m\right)$ we find that the right hand side equals
\begin{align*}
\sum_{\substack{r|\nu \\ r|a}}\epsilon_p(r)\frac{r}{a} \sum_{d\,(a/r)^*} e\left( \frac{dp\nu\nu'/r^2}{a/r}\right)
\sum_{m\,(a)}e\left(-m\frac{\bar d r+h}{a}\right).
\end{align*}
The sum over $m\,(a)$ vanishes unless $\bar d r+h\equiv0\pmod{a}$ in which case it equals $a$. But  $\bar d
r+h\equiv0\pmod{a}$ implies that $h\equiv 0\pmod{r}$ and $\bar d +h/r\equiv0\pmod{a/r}$. Since $\bar d$ is coprime
to $a/r$ we find that $h/r$ must also be coprime to $a/r$ and consequently $r=(h,a)$.

Thus, denoting $r=(h,a)$ and by
$\bar h_1$ a multiplicative inverse of $h/r$ modulo $a/r$ we obtain that the right hand side of \eqref{eq:key2}
equals
\begin{align} \label{eq:key3}
\epsilon_p(r)r \cdot e\left( -\frac{\bar h_1 p\nu\nu'/r^2}{a/r}\right).
\end{align}

We now consider the left hand side of \eqref{eq:key2}. We first assume that $(p,a)=1$. In this case it is equal to
\begin{align*}
&\frac{\epsilon_p(a)}{ a }\sum_{m\,(a)}\sum_{\substack{\lambda\in \calO_F/a\calO_F\\ \lambda\lambda'=-mp\,(a)}} e\left(\frac{\tr(\lambda \nu)-hm}{a}\right) \\
&=\frac{\epsilon_p(a)}{ a }\sum_{\substack{\lambda\in \calO_F/a\calO_F}} e\left(\frac{\tr(\lambda \nu)+h\bar p
\lambda\lambda'}{a}\right).
\end{align*}
Here $\bar p$ denotes a multiplicative inverse of $p$ modulo $a$. We first observe that the sum vanishes unless
$r=(h,a)$ divides $\nu$. This is easily seen by replacing $\lambda\mapsto \lambda+\frac{a}{r} \tau$ for $\tau\in
\calO_F$ in the sum. Therefore, put $r=(h,a)$ and write $h=h_1 r$, $a= a_1 r$, and $\nu=\nu_1 r$. Then the sum
becomes
\begin{align*}
&\frac{\epsilon_p(a)}{ a }\sum_{\substack{\lambda\in \calO_F/a\calO_F}} e\left(\frac{\tr(\lambda \nu_1)+h_1\bar p \lambda\lambda'}{a_1}\right)\\
&=\frac{\epsilon_p(a)}{ a }r^2 e\left(-\frac{p \bar h_1 \nu_1\nu_1'}{a_1}\right)
\sum_{\substack{\lambda\in \calO_F/a_1\calO_F}} e\left(\frac{h_1\bar p(\lambda +p\bar h_1 \nu_1')(\lambda' +p\bar h_1 \nu_1)}{a_1}\right)\\
&=\frac{\epsilon_p(a)}{ a_1 }r e\left(-\frac{p \bar h_1 \nu_1\nu_1'}{a_1}\right)\sum_{\substack{\lambda\in
\calO_F/a_1\calO_F}} e\left(\frac{h_1\bar p\lambda \lambda' }{a_1}\right).
\end{align*}
The latter sum is computed in \cite{Za} \S4 Lemma 2. Inserting its value
\begin{align*}
\frac{1}{ a_1 }\sum_{\substack{\lambda\in \calO_F/a_1\calO_F}} e\left(\frac{h_1\bar p\lambda \lambda'
}{a_1}\right)= \epsilon_p( a_1)
\end{align*}
we finally find that the  left hand side of \eqref{eq:key2} is equal to
\begin{align*}
\epsilon_p(r) r \cdot e\left(-\frac{p \bar h_1 \nu\nu'/r^2}{a/r}\right).
\end{align*}

Let us now consider the left hand side of \eqref{eq:key2} in the case that $p\mid a$. Then it is given by
\begin{align*}
&\frac{1}{ a }\sum_{m\,(a)}\sum_{\substack{\lambda\in \calO_F/a\calO_F\\ \lambda\lambda'\equiv-mp\,(a)}}
\epsilon_p\left(\frac{\lambda\lambda'+mp}{a}\right)  e\left(\frac{\tr(\nu\lambda)-mh}{ a}\right).
\end{align*}
We substitute  $m=m_1+m_2a/p$, where  $m_1$ runs modulo $a/p$ and $m_2$ modulo $p$. Moreover, we notice that in
the sum over $\lambda$, we actually only sum over $\lambda\in \partial_F/a\calO_F$, since $p\mid a$. Substituting
$\lambda\mapsto \sqrt{p}\lambda$, we obtain
\begin{align*}
&\frac{1}{ a }\sum_{\substack{m_1\,(a/p)\\m_2\,(p)}}
\sum_{\substack{\lambda\in \calO_F/a\partial_F^{-1}\\ \lambda\lambda'\equiv m_1\,(a/p)}} \epsilon_p\left(\frac{-\lambda\lambda'+m_1}{a/p}+m_2\right)  e\left(\frac{\tr(\sqrt{p}\lambda\nu)-m_1 h}{ a}\right)e\left(-\frac{hm_2}{p}\right)\\
&=\frac{1}{ ap }\sum_{\substack{m_2\,(p)}}
\sum_{\substack{\lambda\in \calO_F/a\calO_F}}  e\left(\frac{\tr(\sqrt{p}\lambda\nu)-h\lambda\lambda'}{ a}\right)\epsilon_p\left(m_2\right) e\left(-\frac{hm_2}{p}\right)\\
&=\frac{1}{ a\sqrt{p} }\epsilon_p(h) \sum_{\substack{\lambda\in \calO_F/a\calO_F}}
e\left(\frac{\tr(\sqrt{p}\lambda\nu)-h\lambda\lambda'}{ a}\right) .
\end{align*}
Here, in the last line we have inserted the value of the Gauss sum. In particular we see that the latter quantity
vanishes if $p\mid h$. If $p\nmid(h,a)$ we see by replacing  $\lambda\mapsto \lambda+\frac{a}{r} \tau$ for
$\tau\in \calO_F$ that the latter quantity actually vanishes unless $(h,a)$ divides $\nu$. Therefore, as above, we
put $r=(h,a)$ and write $h=h_1 r$, $a= a_1 r$, and $\nu=\nu_1 r$. Then the sum becomes
\begin{align*}
&\frac{r^2}{ a\sqrt{p} }\epsilon_p(h)
\sum_{\substack{\lambda\in \calO_F/a_1\calO_F}}  e\left(\frac{\tr(\sqrt{p}\lambda\nu_1)- h_1\lambda\lambda'}{ a_1}\right)\\
&=\frac{r}{ a_1\sqrt{p} }\epsilon_p(h)e\left(-\frac{\bar h_1 p \nu_1\nu_1'}{a_1}\right)
\sum_{\substack{\lambda\in \calO_F/a_1\calO_F}}  e\left(-\frac{h_1(\lambda-(\bar h_1\sqrt{p}\nu_1)')(\lambda'-(\bar h_1\sqrt{p}\nu_1)')'}{ a_1}\right)\\
&=\frac{r}{ a_1\sqrt{p} }\epsilon_p(h)e\left(-\frac{\bar h_1 p \nu\nu'/r^2}{a/r}\right) \sum_{\substack{\lambda\in
\calO_F/a_1\calO_F}}  e\left(-\frac{h_1\lambda\lambda'}{ a_1}\right).
\end{align*}
By \cite{Za} \S4 Lemma 2 (noting that $p\nmid r$ implies $p\mid a_1$) we have
\begin{align*}
\frac{1}{ a_1 }\sum_{\substack{\lambda\in \calO_F/a_1\calO_F}} e\left(-\frac{h_1\lambda \lambda' }{a_1}\right) =
\epsilon_p( h_1)\sqrt{p}.
\end{align*}
We finally find that left hand side of \eqref{eq:key2} is equal to
\begin{align*}
\epsilon_p(r) r \cdot e\left(-\frac{p \bar h_1 \nu\nu'/r^2}{a/r}\right).
\end{align*}
This concludes the proof of the Lemma.
\end{proof}

We are now ready to compute the Fourier expansion of $\tilde\Phi_{mp^2}(z_1,z_2,s)$ in terms of the coefficients
of the Poincar\'e series $F_m(\tau,s)$.

\begin{theorem}\label{fourier1}
The automorphic Green function $\tilde\Phi_{mp^2}(z_1,z_2,s)$ associated to $\tilde T_m$ has the Fourier expansion
\begin{align*}
&\tilde\Phi_{mp^2}(z_1,z_2,s)\\
&=\left(\frac{p}{\pi}\right)^{s-1/2}\Gamma(s-1/2)L(2s-1,\epsilon_p)
b_m(0,s)(y_1 y_2)^{1-s}\\
&\phantom{=}{}+\sqrt{p} \sum_{\substack{\nu\in \partial_F^{-1}\\\nu\nu'\neq 0}}\sum_{n\geq 1}
 \frac{\epsilon_p(n)}{n} b_m(p\nu\nu',s)\calK_{s}(2\pi n\nu y_1) \calK_{s}(2\pi n\nu' y_2)e(\nu n x_1+\nu'n x_2)\\
&\phantom{=}{}+\sqrt{p} \sum_{\substack{ \lambda\in\partial_F^{-1}\\ \norm(\lambda)=-m/p}} \sum_{n\geq 1}
\frac{\epsilon_p(n)}{n}\big(
2\calI_{s}(2\pi n \beta(\lambda y_1,\lambda' y_2))+\calK_{s}(2\pi n \beta(\lambda y_1,\lambda' y_2))\big)\\
&\qquad\qquad\qquad\qquad\qquad\qquad\qquad \times \calK_{s}(2\pi n \alpha(\lambda y_1,\lambda' y_2)) e(n\lambda
x_1+n\lambda' x_2).
\end{align*}
If converges normally for $y_1 y_2>mp$. Here $b_m(n,s)$ denote the Fourier coefficients of the Poincar\'e series $F_m(\tau,s)$, and  $L(s,\epsilon_p)$ the
Dirichlet $L$-function corresponding to $\epsilon_p$.
\end{theorem}

\begin{proof}
We use \eqref{Phim2}, Lemma \ref{lphi0}, and Lemma \ref{lHA} to compute the Fourier expansion. The constant term
is given by
\begin{align*}
& 2  \sum_{a=1}^{\infty} \tilde G_a(mp^2,0)b_s^{mp/a^2}(0,y_1,y_2)\\
&=\frac{\pi \Gamma(s-1/2)^2}{ \sqrt{p}\,\Gamma(2s)}(4mp)^s (y_1 y_2)^{1-s} \sum_{a=1}^{\infty} \tilde G_a(mp^2,0)
a^{-2s}.
\end{align*}
In view of Lemma \ref{GaHb} this is equal to
\begin{align}
\nonumber &\frac{\pi \Gamma(s-1/2)^2}{ \sqrt{p}\,\Gamma(2s)}(4mp)^s (y_1 y_2)^{1-s} \sum_{a=1}^{\infty} \sum_{r|a}
\epsilon_p(r) H_{a/r}\left( m,0\right)
 a^{1-2s}\\
\nonumber &=\frac{\pi \Gamma(s-1/2)^2}{ \sqrt{p}\,\Gamma(2s)}(4mp)^s (y_1 y_2)^{1-s} \sum_{c=1}^{\infty}
\sum_{r=1}^\infty \epsilon_p(r) H_{c}\left( m,0 \right)
 (cr)^{1-2s}\\
\label{h1} &=\frac{\pi \Gamma(s-1/2)^2}{ \sqrt{p}\,\Gamma(2s)}(4mp)^s (y_1 y_2)^{1-s} L(2s-1,\epsilon_p)
\sum_{c=1}^{\infty} H_{c}\left( m,0\right) c^{1-2s}.
\end{align}
By means of the formula for the constant coefficient $b_m(0,s)$ of $F_m(\tau,s)$ and the duplication formula
$\Gamma(s)\Gamma(s+\frac{1}{2})=2^{1-2s}\sqrt{\pi}\Gamma(2s)$ we see that the constant term equals:
\begin{align*}
\left(\frac{p}{\pi}\right)^{s-1/2}\Gamma(s-1/2)L(2s-1,\epsilon_p) b_m(0,s)(y_1 y_2)^{1-s}.
\end{align*}


We now consider the $\nu$-th Fourier coefficient of $\tilde\Phi_{mp^2}(z_1,z_2,s)$ for $\nu\in \partial_F$ with
$\nu\nu'>0$. It is given by
\begin{align*}
& 2  \sum_{a=1}^{\infty} \tilde G_a(mp^2,\nu)b_s^{mp/a^2}(\nu,y_1,y_2)\\
&=2\pi  \sum_{a=1}^{\infty} \frac{1}{a}\tilde
G_a(mp^2,\nu)\sqrt{\frac{m}{\nu\nu'}}\,I_{2s-1}\left(\frac{4\pi}{a}\sqrt{mp|\nu\nu'|}\right)\calK_{s}(2\pi \nu
y_1) \calK_{s}(2\pi \nu' y_2).
\end{align*}
In view of Lemma \ref{GaHb} this is equal to:
\begin{align*}
&2\pi  \sum_{a=1}^{\infty}  \sum_{\substack{r|\nu \\ r|a}}
\epsilon_p(r) H_{a/r}\left(m, \frac{p\nu\nu'}{ r^2}\right)\sqrt{\frac{m}{\nu\nu'}}\,I_{2s-1}\left(\frac{4\pi}{a}\sqrt{mp|\nu\nu'|}\right)\calK_{s}(2\pi \nu y_1) \calK_{s}(2\pi \nu' y_2)\\
&=2\pi  \sqrt{p}  \sum_{\substack{r|\nu }} \frac{\epsilon_p(r)}{r}
\sqrt{\frac{mr^2}{p\nu\nu'}}\sum_{c=1}^{\infty}   H_{c}\left(m, \frac{p\nu\nu'}{ r^2}\right)
I_{2s-1}\left(\frac{4\pi}{cr}\sqrt{mp|\nu\nu'|}\right)\calK_{s}(2\pi \nu y_1) \calK_{s}(2\pi \nu' y_2).
\end{align*}
By Theorem \ref{fourierF} we finally find for the $\nu$-th coefficient:
\begin{align*}
&  \sqrt{p}  \sum_{\substack{r|\nu }} \frac{\epsilon_p(r)}{r} b_m(p\nu\nu'/r^2,s)\calK_{s}(2\pi \nu y_1)
\calK_{s}(2\pi \nu' y_2).
\end{align*}

In the same way one can show that for $\nu\in \partial_F^{-1}$ with $\nu\nu'<0$ the $\nu$-th coefficient of
$\tilde\Phi_{mp^2}(z_1,z_2,s)-\tilde\Phi_{mp^2}^0(z_1,z_2,s)$ is equal to
\begin{align*}
&  \sqrt{p}  \sum_{\substack{r|\nu }} \frac{\epsilon_p(r)}{r} b_m(p\nu\nu'/r^2,s)\calK_{s}(2\pi \nu y_1) \calK_{s}(2\pi \nu' y_2)\\
&+ \sqrt{p}  \sum_{\substack{r|\nu }} \frac{\epsilon_p(r)}{r} \delta_{-m,p\nu\nu'/r^2}\calK_{s}(2\pi \nu y_1)
\calK_{s}(2\pi \nu' y_2).
\end{align*}
Here the extra contribution with the Kronecker delta comes from the Kronecker delta $\delta_{-m,n}$ in the formula
for $b_m(n,s)$ with $n<0$.

If we put the above contributions together, and use in addition the formula for $\tilde\Phi_{mp^2}^0$ of Lemma
\ref{lphi0}, we find
\begin{align*}
&\tilde\Phi_{mp^2}(z_1,z_2,s)\\
&= \left(\frac{p}{\pi}\right)^{s-1/2}\Gamma(s-1/2)L(2s-1,\epsilon_p)
b_m(0,s)(y_1 y_2)^{1-s}\\
&\phantom{=}{}+\sqrt{p} \sum_{\substack{\nu\in \partial_F^{-1}\\\nu\nu'\neq 0}}
\sum_{\substack{r|\nu }} \frac{\epsilon_p(r)}{r} b_m(p\nu\nu'/r^2,s)\calK_{s}(2\pi \nu y_1) \calK_{s}(2\pi \nu' y_2)e(\nu x_1+\nu' x_2)\\
&\phantom{=}{}+\sqrt{p}  \sum_{\substack{\nu\in \partial_F^{-1}\\\nu\nu'< 0}} \sum_{\substack{r|\nu }}
\frac{\epsilon_p(r)}{r} \delta_{-m,p\nu\nu'/r^2}\calK_{s}(2\pi \nu y_1)
\calK_{s}(2\pi \nu' y_2)e(\nu x_1+\nu' x_2)\\
&\phantom{=}{}+2\sqrt{p} \sum_{\substack{ \lambda\in\partial_F^{-1}\\ \norm(\lambda)=-m/p}} \sum_{n\geq 1}
\frac{\epsilon_p(n)}{n}
 \calI_{s}(2\pi n \beta(\lambda y_1,\lambda' y_2))
\calK_{s}(2\pi n \alpha(\lambda y_1,\lambda' y_2)) e(n\lambda x_1+n\lambda' x_2).
\end{align*}
By rearranging the sums one deduces the stated formula.
\end{proof}

Now it can be proved as in \cite{Br1} that
$\tilde\Phi_{mp^2}(z_1,z_2,s)$ has a meromorphic continuation in
$s$ to a neighborhood of $s=1$. The continuation turns out to be
holomorphic at $s=1$, whereas there always was a simple pole in
\cite{Br1}. This follows from the presence of $L(2s-1,\epsilon_p)$
in the constant term of the Fourier expansion given in Theorem
\ref{fourier1}, while there appeared $\zeta(2s-1)$ in \cite{Br1}
Theorem 1.

\begin{definition}
We define the regularized Green function $\tilde\Phi_{mp^2}(z_1,z_2)$ for $\tilde T_m$ as the value of
$\tilde\Phi_{mp^2}(z_1,z_2,s)$ at $s=1$.
\end{definition}

One finds that $\tilde\Phi_{mp^2}(z_1,z_2)$ is a harmonic function on $\H^2\setminus\tilde T_m$ with a logarithmic
singularity along $-2\tilde T_m$.

\begin{theorem}\label{fourier2}
The Green function $\tilde\Phi_{mp^2}(z_1,z_2)$ associated to $\tilde T_m$ has the Fourier expansion
\begin{align*}
&\tilde\Phi_{mp^2}(z_1,z_2)\\
&=\sqrt{p}L(1,\epsilon_p)
b_m(0,1)\\
&\phantom{=}{}+\sqrt{p} \sum_{\substack{\nu\in \partial_F^{-1}\\\nu\nu'> 0}}\sum_{n\geq 1}
 \frac{\epsilon_p(n)}{n} b_m(p\nu\nu',1)e(-2\pi n|\nu y_1+\nu' y_2|)e(\nu n x_1+\nu'n x_2)\\
&\phantom{=}{}+\sqrt{p} \sum_{\substack{ \lambda\in\partial_F^{-1}\\ \norm(\lambda)=-m/p}} \sum_{n\geq 1}
\frac{\epsilon_p(n)}{n} e\big(-2\pi n |\lambda y_1+\lambda' y_2|\big) e(n\lambda x_1+n\lambda' x_2).
\end{align*}
It converges normally on $y_1 y_2>mp$ outside the polar part of $\tilde T_m$.
\end{theorem}

\begin{proof}
This follows immediately from Theorem \ref{fourier1}, \eqref{Mspecial}, and \eqref{Wspecial}, noting that
$b_m(n,1)=0$ when $n<0$, and that $\alpha(r_1,r_2)-\beta(r_1,r_2)=|r_1+r_2|$ for $r_1,r_2\in \R$ with $r_1r_2<0$.
\end{proof}

\begin{remark}
We have $\sqrt{p}L(1,\epsilon_p)=h_F \log(\eps_0)$, where $h_F$
denotes the class number of $F$ and $\eps_0>1$ the fundamental unit.
\end{remark}

\section{Twisted Borcherds products}
\label{sect:5}

As an application of Theorem \ref{fourier2} of the previous section we
obtain a variant of the Borcherds lift for Hilbert modular
surfaces (see \cite{Bo2}, \cite{Br1}, \cite{BB}). It can be viewed as a
multiplicative analogue of the Doi-Naganuma lift \cite{DN} from
holomorphic modular forms of weight $k$ for $\Gamma'$ to Hilbert
modular forms of weight $k$ for $\Gamma$. Its existence was suggested
by Zagier in \cite{Za3}.

Following \cite{Za3} \S7, we define
\begin{align}\label{defRp}
R_p(t)=\prod_{b\,(p)} (1-e(b/p)t)^{\epsilon_p(b)}.
\end{align}
It is a rational function of $t$ with coefficients in $F$.

\begin{lemma}\label{Rproperties}
The function $R_p(t)\in F(t)$ has the following properties:
\begin{enumerate}
\item[(i)] $R_p(t)'=R_p(t)^{-1}$.
\item[(ii)] $R_p(t^{-1})=R_p(t)$.
\item[(iii)] $\log (R_p(t))= -\sqrt{p}\sum_{n\geq 1}\frac{1}{n}\epsilon_p(n)t^n$.
\end{enumerate}
\end{lemma}

\begin{proof} The first two properties are verified by direct computation. The third property follows from the identity $\sum_{b\,(p)}\epsilon_p(b) e(bn/p)=\sqrt{p}\epsilon_p(n)$.
\end{proof}

\begin{theorem}\label{hilbert}
Let $f=\sum_{n\gg -\infty}c(n)q^n\in \Z[j]$ be a weakly holomorphic
modular form of weight $0$ for $\Gamma'$ with integral Fourier
coefficients. Then there exists a symmetric meromorphic Hilbert
modular function $\Psi(z,f)$ for $\Gamma$ (of weight $0$, with trivial
multiplier system, defined over $F$) such that:
\begin{enumerate}
\item[(i)] The divisor of $\Psi(z,f)$ is given by
\[ \dv(\Psi(z,f))= \sum_{n>0}  c(-n) \tilde T_{n}.\]
\item[(ii)]
The function $\Psi(z,f)$ has the Borcherds product expansion
\begin{align*}
\Psi(z,f)&= \prod_{\substack{\nu\in\partial_F^{-1} \\
\nu>0}}\prod_{b\,(p)} \left(1-e(\tfrac{b}{p}+\nu
z_1+\nu'z_2)\right)^{\epsilon_p(b)c(p\nu\nu')}\\
&=\prod_{\substack{\nu\in\partial_F^{-1} \\ \nu>0}} R_p\big(e(\nu
z_1+\nu'z_2)\big)^{c(p\nu\nu')} ,
\end{align*}
which converges normally for all $(z_1,z_2)$ with $y_1 y_2 > Np$
outside the set of poles, where $N=\max \{n\in \Z;\; c(-n)\neq 0\}$.
\item[(iii)] The lifting is multiplicative, i.e., if $f,g\in \Z[j]$,
then $\Psi(f+g)=\Psi(f)\Psi(g)$.
\item[(iv)]
We have
\[
\log|\Psi(z,f)| =-\frac{1}{2} \sum_{n>0} c(-n) \left(
\tilde\Phi_{np^2}(z)-\sqrt{p} L(1,\epsilon_p) b_{n}(0,1)\right).
\]
\end{enumerate}
\end{theorem}

\begin{proof}
Let us first assume that $f=q^{-m}+O(1)$, as $q\to 0$, for some positive
integer $m$.  Then, according to Theorem \ref{fourierF1},
$f(\tau)=F_m(\tau,1)+C$ for some constant $C\in \Z$.  We define the
function $\Psi(z,f)$ by the product expansion in (ii).  In view of
Theorem~\ref{fourier2} and Lemma~\ref{Rproperties} (ii), (iii) we have
\begin{align}\label{intid}
\log|\Psi(z,f)| =-\frac{1}{2} \left( \tilde\Phi_{mp^2}(z)-\sqrt{p}
L(1,\epsilon_p) b_{m}(0,1)\right).
\end{align}
In particular, the normal convergence of the Fourier expansion of
$\tilde\Phi_{mp^2}(z)$ on $y_1 y_2>mp$ implies the normal convergence
of the infinite product $ \Psi(z,f)$. Therefore (ii) and (iv) hold.

In the same way as in \cite{Br1}, Theorem 4, one can show that $
\Psi(z,f)$ has a meromorphic continuation to all of $\H^2$ and that
$\dv(\Psi(z,f))=\tilde T_m$: This follows from \eqref{intid} and the
fact that $\tilde\Phi_{mp^2}(z)$ is a pluriharmonic function on
$\H^2\setminus \tilde T_m$ with a logarithmic singularity along
the divisor $-2\tilde T_m$.

Moreover, \eqref{intid} implies that $|\Psi(z,f)|$ is invariant under
the group $\Gamma$. The matrix $S=\kzxz{0}{-1}{1}{0}\in \Gamma$
satisfies the relation $S^2=-1$. Consequently,
\begin{align}\label{intid2}
\Psi(Sz,f)=\pm \Psi(z,f).
\end{align}
On the other hand, it follows from the product expansion that
$\Psi(z,f)$ is invariant under translations $\kzxz{1}{\mu}{0}{1}$
where $\mu\in \OK$.  Using the relation $(S\kzxz{1}{1}{0}{1})^3=-1$ in
$\Gamma$, we may conclude that the sign in \eqref{intid2} must be
$+1$.  But, by a theorem of Vaserstein, the translations
$\kzxz{1}{\mu}{0}{1}$ and the matrix $S=\kzxz{0}{-1}{1}{0}$ generate the group
$\Gamma$. Hence $\Psi(z,f)$ is invariant under $\Gamma$.

Moreover, the product expansion implies that $\Psi(z,f)$ can be
written as the quotient of two holomorphic Hilbert modular forms with
Fourier coefficients in $F$. Hence, by the $q$-expansion principle, $\Psi(z,f)$
is defined over $F$.

Finally, the multiplicativity (iii) also follows from the infinite
product expansion.
%
%
\end{proof}

The condition $\nu>0$ under the product replaces the Weyl-chamber
condition occurring for the untwisted Borcherds products \cite{BB}.
It has this easy form because of the presence of the character
$\epsilon_p(b)$ which causes some cancellations. Because of Lemma \ref{Rproperties} (ii) the product is actually over $\calO_F^\times/\{\pm 1\}$.

For the rest of this paper we write $\tilde\Psi_m$ ($m\in\Z_{>0}$)
for the  twisted Borcherds
lift in the sense of Theorem \ref{hilbert} of the unique weakly
holomorphic modular form
\[
J_m=q^{-m}+O(q)\in \Z[j].
\]

\begin{remark}\label{rem:conj}
In view of Lemma \ref{Rproperties} (i), the conjugation in $F$ maps $\Psi(z,f)$ to $\Psi(z,f)^{-1}$.
\end{remark}

\section{Class Fields}
\label{sect:class}

Let $F=\mathbb Q(\sqrt p)$ be a real quadratic field as before.
Let $K = F(\sqrt\Delta)$ be a non-biquadratic totally imaginary
quadratic extension of $F$.  We view both $K$ and
$F(\sqrt{\Delta'})$ as subfields of $\mathbb C$ with
$\sqrt{\Delta}, \sqrt{\Delta'} \in \mathbb H$. Then
$M=F(\sqrt\Delta, \sqrt{\Delta'})$ is Galois over $\Q$ and has an
automorphism $\sigma$ of order $4$ such that
$\sigma(\sqrt{\Delta}) =\sqrt{\Delta'}$ and
$\sigma(\sqrt{\Delta'})=-\sqrt\Delta$.  The field $K$ has four CM
types: $\Phi=\{1, \sigma \}$, $\sigma\Phi =\{\sigma, \sigma^2\}$,
$\sigma^2\Phi$, and $\sigma^3 \Phi$.  We assume that the relative
discriminant $d_{K/F}$ of $K/F$ satisfies the technical condition
\begin{equation} \label{oldeq1.1}
d_{K/F} \cap \mathbb Z = q \mathbb Z, \quad \norm_{F/\mathbb Q} d_{K/F} =q,
\end{equation}
for a prime number $ q \equiv 1 \pmod 4$.

Recall that the Hilbert modular surface $X$ corresponding to $\Gamma=\Sl_2(\O_F)$
parameterizes isomorphism
classes of triples $(A, \imath, m)$, where $(A, \imath)$ is an abelian
surface with real multiplication $ \imath: \O_F \hookrightarrow
\hbox{End}(A), $ and
$$
m:(\mathfrak M_A, \mathfrak M_A^+) \longrightarrow \left(
\partial_F^{-1}, \partial_F^{-1, +}\right)
$$
is an $\O_F$-isomorphism between the polarization module $\mathfrak
M_A$ of $A$ and $\partial_F ^{-1}$, taking the subset of polarizations
to totally positive elements of $\partial_F^{-1}$ (see e.g.~\cite{Go},
Theorem 2.17 and \cite{BY} Section 3).

Let $\CM(K, \Phi, \calO_F)$ be the CM $0$-cycle in  $X$ of CM abelian surfaces
of CM type $(K, \Phi)$, i.e., the points on $X$ with an $\calO_K$-action  via $\Phi$ (see \cite{BY} Section 3 for details).
The field of moduli for the CM cycle
$\CM(K):=\CM(K, \Phi,
\calO_F) + \CM(K, \sigma^3\Phi, \calO_F)$ is $\mathbb Q$.

Let $(\tilde K, \tilde\Phi)$ be the reflex of $(K, \Phi)$ with maximal
totally real subfield $\tilde F =\mathbb Q(\sqrt{\Delta \Delta'})$.
Let $I(K)$ be the group of all fractional ideals of $K$, and let $H(
K)$ be the subgroup of $I( K)$ of ideals $ {\mathfrak a}$ such that
    \begin{equation}
    N_{  \Phi}  {\mathfrak a} = \mu \O_{\tilde K},
   \quad N  {\mathfrak a} =\mu \bar\mu \quad \text{ for some $\mu \in  \tilde K^*$}.
   \end{equation}
  Here $N\mathfrak a =\# \O_K/\mathfrak a$, and $N_{ \Phi}$ is the type norm
  from $I(  K)$ to $I(\tilde K)$ given by, in this case,
  $$
  N_{  \Phi} (\mathfrak a) = N_{M/\tilde K} (\mathfrak a \O_M).
  $$
   We call the quotient $\CC(K)= \CC(  K, \Phi) = I(  K)/H(  K)$
   the CM ideal class group of $  K$. According to
   \cite{ShimuraCM}, page 112, Main Theorem 1, the class field
   $H_{\tilde K}$
   of $\tilde K$ associated to the CM ideal class group
   $\CC(\tilde K,\tilde\Phi)$ is the composite of $\tilde K$ with the field
   of the moduli of any polarized CM abelian variety of type
   $(K, \Phi)$ by $\O_K$.
Recall   that
   $\CC(\tilde K, \tilde\Phi) =\Gal(H_{\tilde K}/\tilde K)$ acts
   on $\CM(K, \Phi, \O_F)$ via \cite{BY}, (3.7).

Inspired by the classical result that
   $j(\frac{-D +\sqrt{-D}}2)$ generates the Hilbert class field of
   the imaginary quadratic field $\mathbb Q(\sqrt{-D})$, we
   consider whether $H_{\tilde K}$ has a `canonical' generator
   over $\tilde K$.  Since $F$ is in general not contained in
   $\tilde K$, it seems more natural to consider the field
   $L_K= K H_{\tilde K} = M H_{\tilde K} =F H_{\tilde K}$ over
   either $M$ or $F$, where $M=K \tilde K$ is the smallest Galois
   extension of $\mathbb Q$ containing $K$ and/or $\tilde K$. For
   example, a natural question is whether $L_K=M(\tilde\Psi_m(z))$
   for some $m$ and some CM point $z \in \CM(K, \Phi, \O_F)$?


\begin{lemma} \label{cllem1} Assume that $d_K=p^2 q$ with $p \equiv q \equiv 1
\mod 4$ being odd primes. Then $\CC(\tilde K, \tilde\Phi) $ acts
on $\CM(K, \Phi, \O_F)$ simply transitively. In particular,
$\CM(K, \Phi, \O_F)$ is a single Galois orbit of a CM point.
\end{lemma}

\begin{proof}
Let $\CL_0(K)$ be the subgroup of $\CL(K)$ generated by
ideal classes  $[\mathfrak a]$ such that
$$
N_{K/F} \mathfrak a  =\mu \O_F, \quad \text{$\mu \gg 0$ (totally positive)}.
$$
Then \cite{BY} Lemma 5.3 asserts that the type norm  $N_{\tilde
\Phi}$ gives an isomorphism between  $\CL_0(\tilde K)$ and
$\CL_0(K)$. It is clear by definition that $N_{\tilde \Phi}$ maps
$I(\tilde K)$ to $\CL_0(K)$, so it is surjective. It is easy to
check that the kernel is exactly $H(\tilde K)$. Indeed, if
$[\mathfrak a] \in \ker N_{\tilde \Phi}$, then
$$
N_{\tilde \Phi} \mathfrak a  =\mu \O_K.
$$
This implies
$$
N \mathfrak a  \O_F =N_{K/F} N_{\tilde \Phi} \mathfrak a= \mu \bar
\mu \O_F,
$$
and thus
$$
N \mathfrak a = \mu \bar \mu \epsilon
$$
for some unit $\epsilon$. Clearly $\epsilon$ is totally positive
and thus a square $\epsilon_1^2$ since $F=\mathbb Q(\sqrt p)$.
Replacing $\mu$ by $\mu \epsilon_1$, one sees that $\mathfrak a
\in H(\tilde K)$. So
$$
\CC(\tilde K) \cong \CL_0(K),
$$
and $[L_K: M] = [H_{\tilde K}: \tilde K] =\# \CL_0(K)$. On the
other hand, according to   the remark after  \cite{BY}, Lemma 3.3,
The forgetful map
$$\CM(K, \Phi, \calO_F) \rightarrow \CL_0(K), \quad
(\mathbb C^2/\Phi(\mathfrak a), \imath, m) \rightarrow \mathfrak a$$ is a bijection.   So $\CC(\tilde K)$ acts on
$\CM(K, \Phi, \O_F)$ simply transitively.
\end{proof}

 \begin{corollary}
 Let the assumption be as in Lemma~\ref{cllem1}.
 Let $z =(A, \imath, m) \in \CM(K, \Phi, \O_F)$ and let $k_z$ be the
 field of definition of $z$ containing $\tilde F$.  Then
 $$
 \tilde\Psi_m (\CM(K, \Phi, \O_F))
   = N_{M_z/M^+} \tilde\Psi_m(z).
   $$
   Here $M_z= F k_z=M^+ k_z$, and $M^+ =F \tilde F$ is the maximal
   totally real subfield of $M$.
   \end{corollary}

\begin{proof}
Recall that $k_z$ is the field of moduli of $(A, \imath, m)$ and that $H_{\tilde K} =\tilde K k_z$. Lemma~3.4(1) of
\cite{BY} asserts that $\CM(K, \Phi, \O_F)$ is defined over $\tilde F$. So the lemma above implies that $k_z$ does
not contain $\tilde K$ and $[k_z: \tilde F] =\#\CC(\tilde K)$, and different embeddings of $k_z$ into $\mathbb C$
fixing $\tilde F$ map $z$ into different CM points in $\CM(K, \Phi, \O_F)$. Since $\tilde\Psi_m$ is defined over
$F$, we obtain the corollary.
\end{proof}

\begin{corollary}
\label{clcor1.3}
Let the assumption be as in  Lemma~\ref{cllem1}.
 If $\tilde\Psi_m(\CM(K))$ is not a $k$-power in $F$
for any odd integer $k >1$, then $\tilde\Psi_m(z)$ generates $L_K$ over $M$.
\end{corollary}

\begin{proof}
Recall that $\CM(K) =\CM(K, \Phi, \O_F) + \CM(K, \sigma^3\Phi, \O_F)$. If $\tilde\Psi_m(z)$ does not generate
$L_K$ over $M$, then it does not generate $M^+ k_z$ over $M^+$. Let $L$ be the subfield of $M^+ k_z$ generated by
$\tilde\Psi_m(z)$ over $M^+$.  Then  the above corollary implies that
$$
\tilde\Psi_m(\CM(K, \Phi, \O_F)) = \left( N_{L/M^+}(\tilde\Psi_m(z))\right)^{[M^+k_z: L]}
$$
is a $k$-th power with $k=[M^+k_z: L]>1$. Therefore
$$
\tilde\Psi_m(\CM(K)) =N_{M^+/F}\tilde\Psi_m(\CM(K, \Phi, \O_F))
   =\left(N_{L/F}(\tilde\Psi_m(z))\right)^{k}.
   $$
   Finally, the condition  $d_K=p^2 q$ implies that $d_{\tilde
   K} =q^2 p$ is an odd number and thus  its class number $h({\tilde K})$
   is odd \cite{ClassNumberParity}. In
   particular, $k | \#\CC(\tilde K) | h(\tilde K)$ is odd.
   \end{proof}

\begin{remark}  \label{rem:square}
Let the notation be as above, and let
$c= \tilde\Psi_m(\CM(K, \Phi, \O_F))  \in M^+$. The above proof implies
$$
N_{M^+/F} c = \tilde\Psi_m(\CM(K)).
$$
On the other hand,  it follows from Remark \ref{rem:conj} that
$$
N_{M^+ /\tilde F} c = 1.
$$
The numerical examples in Section \ref{sect:examples} suggest that $\tilde\Psi_m(\CM(K))$ is a square and that $c
\in F$, that is,
\begin{equation}
\tilde\Psi_m(\CM(K, \Phi, \O_F)) = \tilde\Psi_m(\CM(K, \sigma^3 \Phi, \O_F)).
\end{equation}
\end{remark}

 We will try to compute $\tilde\Psi_m(\CM(K))$ in the next section.

\begin{proposition}\label{prop:pointsep}
Assume that $d_K=p^2 q$ with $p \equiv q \equiv 1 \mod 4$ being odd primes, and let $B$ be a set of positive
integers. Then   $L_K= M(\tilde\Psi_m(z), m \in B)$ for a $z \in \CM(K, \Phi,\calO_F)$ if and only if the functions
$\tilde\Psi_m$ ($m \in B$)  separate the points in $\CM(K, \Phi, \O_F)$.
\end{proposition}

\begin{proof}
First assume that the $\tilde\Psi_m$ ($m \in B$)
 separate the points
in $\CM(K, \Phi, \O_F)$. Clearly $L_K \supset M(\tilde\Psi_m(z), m \in B)$. Suppose that $\alpha \in \Gal(L_K/M) =
\Gal(H_{\tilde K}/\tilde K)$ fixes the field $M(\tilde\Psi_m(z), m \in B)$. Then
$$
\tilde\Psi_m(\alpha(z)) = \alpha (\tilde\Psi_m(z)) =
\tilde\Psi_m(z)
$$
for every $m \in B$, and so $\alpha(z) =z$. This implies that
$\alpha =1$, and $L_K= M(\tilde\Psi_m(z), m \in B)$.

 Conversely, if there are $z_1 \ne z_2 \in \CM(K, \Phi, \O_F)$ such
 that $\tilde\Psi_m(z_1) = \tilde\Psi_m(z_2)$ for every $m \in B$,
 let $1 \ne \alpha \in \Gal(L_K/M)$ such that $\alpha(z_1) =z_2$, which
 exists by Lemma \ref{cllem1}. Then
 $$
 \alpha(\tilde\Psi_m(z_1)) = \tilde\Psi_m(z_2) =\tilde\Psi_m(z_1)
 $$
 for every $m \in B$ and thus $M(\tilde\Psi_m(z), m \in B)$ is
 fixed by $\alpha$. So $L_K \neq M(\tilde\Psi_m(z), m \in B)$.
\end{proof}

 For two positive integers $m$ and $n$ such that $mn$ is not a
 square, it is known  that $\tilde T_{m}$ and $\tilde T_{n}$ has no
 common  component \cite{HZ}. This implies that the rational map
 \begin{equation}
 \phi_{m, n}: X \rightarrow \mathbb P^1 \times \mathbb P^1, \quad
   z \mapsto (\tilde\Psi_m(z), \tilde\Psi_n(z))
   \end{equation}
   is generically finite, i.e.,
   $
   [\mathbb C(X): \mathbb C(\tilde\Psi_m(z),\tilde\Psi_n(z))]
   $
   is finite.

\begin{theorem}
\label{cltheo6.6} Let $F=\mathbb Q(\sqrt p)$ be a fixed real
quadratic field. Let $K$ be a non-biquadratic CM  quartic field
with maximal totally real subfield $F$ and  discriminant $d_K= p^2
q$ with $p \equiv q \equiv 1 \mod 4$ primes. Let $z \in \CM(K,
\Phi, \O_F)$ and $L_{m, n}(K) = M(\tilde\Psi_m(z),
\tilde\Psi_n(z))$ be the unramified abelian extension of $M$
generated by $\tilde\Psi_m(z)$ and $\tilde\Psi_n(z)$. Then there
is a constant  $d
>0$ depending on $F$, $m$, and $n$, but independent of $K$, such
that
\begin{enumerate}
\item[(i)] $[L_K: L_{m, n}(K)] \le d $.
\item[(ii)]
$ \lim_{q \rightarrow \infty} \frac{\log [L_{m, n}(K): M]}{\log
\sqrt q} =1$.
\end{enumerate}
\end{theorem}

\begin{proof}
Lemma \ref{cllem1} implies that $\Gal(L_K/M)$ acts
on $\CM(K, \Phi, \O_F)$  simply transitively, so
$$
\CM(K, \Phi, \O_F) =\{ \sigma(z);\; \sigma \in \Gal(L_K/M)\}.
$$
Here $z$ is a fixed CM point in $\CM(K, \Phi, \O_F)$. Set
\begin{align*}
A&= \{ z'\in \CM(K, \Phi, \O_F);\; \tilde\Psi_i(z')
=\tilde\Psi_i(z),\, i=m, n\}
 \\
 &=\{\sigma(z);\; \sigma\in \Gal(L_K/M), \tilde\Psi_i(\sigma(z))
=\tilde\Psi_i(z),\, i=m, n\}.
\end{align*}
 Since
$$
\sigma(\tilde\Psi_m(z)) = \tilde\Psi_m(\sigma(z)),
$$
one has that $\sigma(z) \in A$ if and only is $\sigma \in
\Gal(L_K/L_{m, n}(K))$. So
$$
[L_K:L_{m, n}(K)] =\# A \le \deg \phi_{m, n}
$$
generically, that is, if $\phi_{m, n}(z) \notin  B= \{ (a, b);\;
\phi_{m, n}^{-1}(a, b) \hbox{ is infinite} \}$. Notice that  $B$
is a finite set. Since $\phi_{m, n}$ is defined over $F$,
$\phi_{m, n}(z) \in B$ for some $z \in \CM(K, \Phi, \O_F)$ implies
that $\phi_{m, n}( \sigma(z)) \in B$ for all $\sigma(z) \in \CM(K,
\Phi, \O_F)$, i.e.,
$$
\CM(K, \Phi, \O_F) \subset \phi_{m, n}^{-1} (B)  \quad \hbox{ or } \quad \CM(K, \Phi, \O_F) \cap \phi_{m, n}^{-1}
(B)=\emptyset. 
$$
According to the equidistribution theorem of CM  points on a
Hilbert modular variety recently proved by Zhang
(cf.~\cite{Zhang}, Theorem 2.1), which extends a well-known
theorem of W. Duke on modular curves \cite{Duke},  $\{ \CM(K,
\Phi, \O_F);\; d_K=p^2 q\}$ is equidistributed on $X$. On the
other hand, $\phi_{m, n}^{-1}(B)$ is just a divisor of $X$, so
there is $q_0
>0$ such that $q > q_0$ implies that $ \CM(K, \Phi, \O_F) \cap
\phi_{m, n}^{-1} (B) $ is empty. Let
$$
d_0 =\max \{[L_K: L_{m, n}(K)];\; q \le q_0\},\quad
\hbox{and}\quad d = \max\{\deg \phi_{m, n}, d_0\}.
$$
Then we always find
$$
[L_K:L_{m, n}(K)] \le d.
$$
This proves (i).

For (ii), one has
$$
[L_K:M] = \#\CM(K, \Phi, \O_F) = \# \CL_0(K) =c \frac{h_K R_K}{h_F
R_F} .
$$
Here $h_F$ and $h_K$ are the class numbers of $F$ and $K$ respectively,  $R_F$ and $R_K$ are the regulators of $F$
and $K$
 respectively, and $\frac{1}2  \le c \le 2$. Indeed,
 $$
 \frac{R_F}{R_K} \le c= [\CL(F): N_{K/F} \CL(K)] \frac{R_F}{R_K}
  \le  [\CL(F): N_{K/F} \CL(K)].
 $$
Since $F$ is fixed, one sees from the Brauer-Siegel theorem
(cf.~\cite{Lang},  Chapter XVI, Lemma~2 and Theorem~5)   that
\begin{equation} \label{cleq6.5}
\lim_{q \rightarrow \infty} \frac{\log [L_K:M]}{\log \sqrt q}=1.
\end{equation}
Now (ii) follows from (i).
\end{proof}

  Bas Edixhoven gave a very nice lower bound for the size of the
  Galois orbit of a CM point on a Hilbert modular surface in general
  (cf.~\cite{Ed}, Section 6).  In particular, (\ref{cleq6.5}) is
  implicitly given in his proof of \cite{Ed}, Theorems 6.2 and 6.4. He
  also pointed out that using a result of Stark (instead of the
  Brauer-Siegel theorem) \cite{Stark}, one gets an effective lower bound $ [L_{m,
  n}(K): M] \gg q^{\frac{1}4}$.

\begin{proposition} The field $L_K$ is Galois over $F$.
\end{proposition}
\begin{proof} By  \cite{ShimuraCM}, page 112, Main Theorem 1,
$H_{\tilde K}$ is the field of moduli of a CM abelian variety $z=(A, \imath, m) \in \CM(K, \Phi, \O_F)$  over
$\tilde K$. So $L_K$ is the field of moduli of a CM abelian variety $z \in \CM(K, \Phi, \O_F)$ over $M$. Let
$\alpha $ be an embedding of $L_K$ into $\bar F$ fixing $F$, and $\alpha|_K$ is either the identity or the complex
conjugation. So $(A^\alpha, i^\alpha, m^\alpha) \in \CM(K, \Phi, \O_F )$ or $\CM(K, \bar\Phi, \O_F)$. Clearly,
$(K, \Phi, \O_F)$ and $(K, \bar\Phi, \O_F)$ have the same reflex field $\tilde K$, and the associated CM ideal
class groups of $\tilde K$ are the same. This implies $\alpha(L_K) = L_K$, i.e.,  $L_K$ is Galois over $F$.
\end{proof}

 In general, $L_K$ is not abelian over $F$. Let $L_F$ be the
 composite of all $L_K$ with $K$ running through non-biquadratic
 CM quadratic extensions of $F$. Then one has an exact sequence
 $$
 1 \rightarrow \mathcal A \rightarrow \Gal(L_F/F) \rightarrow (\mathbb Z/2)^{(\mathbb N)} \rightarrow 0
 $$
 for some abelian group $\mathcal A$. It might  be interesting to study
 the Galois group $\Gal(L_F/F)$. We end this section with the
 following question.

\begin{question}
Is $L_K$ independent of the choice of the CM types of $K$? This is equivalent to the question whether $L_K$ is
Galois over $\mathbb Q$.
\end{question}

\section{Examples}
\label{sect:examples}

It would be very interesting to obtain closed formulas for the values
of twisted Borcherds products at the CM cycles considered in the previous section
in analogy to
\cite{BY}. However, since the twisted Borcherds products are in
general only defined over $F$ (in contrast to the untwisted Borcherds
products, which are essentially defined over $\Q$), their CM values
will lie in $F$. This makes computations
more difficult. Moreover, taking the norm to $\Q$ does not provide any
insight, since $\Psi(z,f)\cdot\Psi(z,f)'=1$ because of Remark
\ref{rem:conj}.
Note that for numeric computations, in the same way as in
\cite{BY}, the problem arises that the product expansion of Theorem
\ref{hilbert} only converges near the cusps. The CM points usually do
not lie in the domain of convergence. Therefore one has to find an
alternative expression for the twisted Borcherds products one wants
to evaluate.

\medskip

Here we discuss some examples in the special case $p=5$ where
$F=\Q(\sqrt{5})$. The fundamental unit of $\calO_F$ is equal to
$\omega=\frac{1+\sqrt{5}}{2}$. The structure of the graded ring of
holomorphic Hilbert modular forms for the group
$\Gamma=\Sl_2(\calO_F)$ was determined by Gundlach \cite{Gu}, see also
\cite{Mu}.  In particular, it turns out that the graded ring
$M_{2*}^{sym}(\Gamma)$ of holomorphic symmetric Hilbert modular forms
of even weight for $\Gamma$ is the polynomial ring
$\C[g_2,g_6,g_{10}]$, where $g_k$ denotes the Eisenstein series (in
the cusp $\infty$ for $\Gamma$) of weight $k$ normalized such that the
constant term is $1$.  Often it is more convenient to replace the
generators $g_6$ and $g_{10}$ by the cusp forms
\begin{align*}
s_6 &= 67\cdot(2^5\cdot 3^3\cdot 5^2)^{-1}\cdot(g_2^3-g_6),\\
s_{10}&= (2^{10}\cdot 3^5\cdot 5^5 \cdot 7)^{-1}\cdot (2^2\cdot 3\cdot 7\cdot 4231 \cdot g_2^5-5\cdot 67\cdot
2293\cdot g_2^2\cdot g_6 +412751\cdot g_{10}).
\end{align*}
We have
\[
M_{2*}^{sym}(\Gamma)=\C[g_2,s_6,s_{10}].
\]
Notice that $g_2,s_6,s_{10}$ all have rational integral and coprime
Fourier coefficients. The cusp form $s_{10}$ is equal to the
(untwisted) Borcherds lift $\Psi_1^2$ in the sense of \cite{BB}.

Recall from Section \ref{sect:5} that $\tilde\Psi_m$ is the symmetric Hilbert modular function of
weight $0$ which is the twisted Borcherds lift in the sense of Theorem
\ref{hilbert} of the unique weakly holomorphic modular form
$J_m=q^{-m}+O(q)\in \Z[j]$.
By Gundlach's theorem, $\tilde\Psi_m$ must be a rational function in
$g_2,s_6,s_{10}$ with coefficients in $F$.
We now discuss how it can be computed.

For simplicity, we assume that $m$ is a {\em square-free} positive
integer.  Then the Hirzebruch-Zagier divisor $T_{mp^2}$ decomposes
into irreducible components $T_{mp^2}=T_m+F^+_{mp^2}+F_{mp^2}^-$.
Recall that $\dv(\tilde \Psi_m)= F_{mp^2}^+-F_{mp^2}^-$. On the other
hand we can construct a symmetric holomorphic Hilbert modular form
with divisor $F_{mp^2}^++F_{mp^2}^-$ of weight $k_m>0$ by taking the
(untwisted) Borcherds lift $H_m=\Psi_{mp^2}/\Psi_m$ in the sense of
\cite{BB} Theorem 9.  For instance, for $m=1$ we have $k_1=60$, and
$H_1$ it is obtained as the lift of the unique weakly holomorphic
modular form $h\in W_0^+(p,\epsilon_p)$ whose Fourier expansion has the
form
\[
h=\frac{1}{2} q^{-25}-q^{-1}+60+438864q+45271325304q^4+\dots.
\]
The product $\tilde \Psi_m \cdot H_m$ is also a symmetric holomorphic
Hilbert modular form of weight $k_m$. Its divisor is equal to $2
F_{mp^2}^+$.  Hence there exist holomorphic Hilbert modular forms
$\Psi_m^+$ and $\Psi_m^-$ of weight $k_m/2$ for $\Gamma$ such that
$(\Psi_m^+)^2=H_m\cdot \tilde \Psi_m$ and $(\Psi_m^-)^2=
H_m/\tilde \Psi_m $.

The function $\Psi_m^\pm$ must also be symmetric, because any skew-symmetric
Hilbert modular form automatically vanishes on $F_1$, contradicting
$\dv(\Psi_m^\pm)=F_{mp^2}^\pm$.
Moreover, $\Psi_m^\pm$ is defined over $F$ and
\begin{align}
(\Psi_m^+)'&=\Psi_m^-,\\
\label{eq:2}
\tilde \Psi_m &=\Psi^+_m/ \Psi_m^-,\\
\label{eq:1}
H_m &= \Psi^+_m\cdot \Psi_m^-.
\end{align}

By Gundlach's theorem, $\Psi_m^\pm$ is a homogeneous polynomial in
$g_2,s_6,s_{10}$ with coefficients in $F$. Using the infinite product
expansions, this polynomial can be determined explicitly.  Its degree
depends on the weight $k_m$.  Unfortunately, it turns out that the $k_m$ are rather
large which makes computations difficult.  The smallest weights
that occur are $k_1=k_2=60$.  For all $m$ the weight $k_m$ is
divisible by $60$, and if $m>2$ then $k_m\geq 120$.

A computation with Maple shows that
\begin{align*}
16 \Psi_1^+ &= (6+2 \sqrt{5}) g_2^{10} s_{10}+(10-2 \sqrt{5}) g_2^9 s_{6}^2+(-2308750-1031750 \sqrt{5}) g_2^7 s_{10} s_{6}\\
&\phantom{=}{}+(-1220450-543450 \sqrt{5}) g_2^6 s_{6}^3+(856853809375+383196837500 \sqrt{5}) g_2^5 s_{10}^2\\
&\phantom{=}{}+(-133751887500-59814018750 \sqrt{5}) g_2^4 s_{10} s_{6}^2\\
&\phantom{=}{}+(-309550426875-138434703750 \sqrt{5}) g_2^3 s_{6}^4\\
&\phantom{=}{}+(23003309053125000+10287392475000000 \sqrt{5}) g_2^2 s_{10}^2 s_{6}\\
&\phantom{=}{}+(-18093694595625000-8091745695000000 \sqrt{5}) g_2 s_{10} s_{6}^3\\
&\phantom{=}{}+(-16048066250700000-7176913583400000 \sqrt{5}) s_{6}^5\\
&\phantom{=}{}+(24527175191718750000+10968886216250000000 \sqrt{5}) s_{10}^3.
\end{align*}
Moreover,
\begin{align*}
\Psi_2^+ &= \scriptstyle
-g_2^{15}+(48072800+21493760 \sqrt{5}) g_2^{12} s_{6}\\
&\phantom{=}{}\scriptstyle
+(-12166677513088000-5441103582617600 \sqrt{5}) g_2^{10} s_{10}\\
&\phantom{=}{}\scriptstyle
+(4809336551424000 \sqrt{5}+10754003449472000) g_2^9 s_{6}^2\\
&\phantom{=}{}\scriptstyle
+(-343213552017810432000000-153489766622216192000000 \sqrt{5}) g_2^7 s_{10} s_{6}\\
&\phantom{=}{}\scriptstyle
+(356570717554737254400000+159463272647655096320000 \sqrt{5}) g_2^6 s_{6}^3\\
&\phantom{=}{}\scriptstyle
+(-1295691015382296818073600000000-579450637646108270264320000000 \sqrt{5}) g_2^5 s_{10}^2\\
&\phantom{=}{}\scriptstyle
+(773852433537584819077120000000 \sqrt{5}+1730386645943677691904000000000) g_2^4 s_{10} s_{6}^2\\
&\phantom{=}{}\scriptstyle
+(-463842956129634468495360000000-207436876158063085617152000000 \sqrt{5}) g_2^3 s_{6}^4\\
&\phantom{=}{}\scriptstyle
+(26896303882903769962250240000000000 \sqrt{5}+60141963825664373337292800000000000) g_2^2 s_{10}^2 s_{6}\\
&\phantom{=}{}\scriptstyle
+(-89933914258012151985733632000000000 \sqrt{5}-201098345763552732433612800000000000) g_2 s_{10} s_{6}^3\\
&\phantom{=}{}\scriptstyle
+(66189203273169170168775966720000000 \sqrt{5}+148003557895357846527777177600000000) s_{6}^5\\
&\phantom{=}{}\scriptstyle
+(-17950573674763054301052928000000000000 \sqrt{5}-40138702971888390552944640000000000000) s_{10}^3.
\end{align*}
These polynomial representations can be used to calculate the Fourier expansions, which in turn can be employed to
compute the values of $\tilde \Psi_1$ and $\tilde \Psi_2$ at CM points.

A pleasant example is the CM point $z_0=(\zeta_5,\zeta_5^2)$ (where $\zeta_5=e^{2\pi i/5}$) corresponding to the
cyclic CM extension $K=\Q(\zeta_5)$ of $F$. It is known that $ z_0$ is an elliptic fixed point of $\Gamma$ of
order $5$. The stabilizer of $z_0$ in $\Gamma$ is the cyclic subgroup generated by $\kzxz{\omega'}{1}{-1}{0}$. This
implies that every Hilbert modular form for $\Gamma$ of weight coprime to $5$ vanishes at $z_0$. In particular
$g_2(z_0)=s_6(z_0)=0$. Consequently, only the term involving $s_{10}^3$ contributes to the value at $z_0$. We find that
\begin{align*}
\tilde \Psi_1(z_0) &=\Psi_1^+(z_0)/\Psi_1^-(z_0) \\
&=\frac{24527175191718750000+10968886216250000000 \sqrt{5}}{24527175191718750000-10968886216250000000 \sqrt{5}}\\
&=\frac{156973921227+70200871784\sqrt{5}}{156973921227-70200871784\sqrt{5}}\\
&=\left(\frac{\omega}{\omega'}\right)^{27}\cdot
\frac{(4+\omega')^5\cdot(5+\omega')^5}{(4+\omega)^5\cdot(5+\omega)^5}.\\
\intertext{Moreover,}
\tilde \Psi_2(z_0) &=\Psi_2^+(z_0)/\Psi_2^-(z_0) \\
&=-\frac{68476004313518731312+30623400094519340937\sqrt{5}}{68476004313518731312-30623400094519340937\sqrt{5}}\\
&=\left(\frac{\omega}{\omega'}\right)^{51}\cdot
\frac{(9+\omega')^5\cdot(10+\omega')^5}{(9+\omega)^5\cdot(10+\omega)^5}.
\end{align*}
Notice that $4+\omega$, $5+\omega$, $9+\omega$, $10+\omega$
are prime elements of $\calO_F$
above $19$, $29$, $89$, $109$, respectively. In particular the same phenomenon
as in \cite{BY} happens: the prime factors of the CM values are small.

\medskip

We would like to evaluate $\tilde \Psi_m$ at other CM cycles as well, say, corresponding to CM extensions $K/F$ such
that $K/\Q$ is non-Galois (and satisfies the assumptions of Section \ref{sect:class}). Here it is rather difficult to compute the value exactly. To get a feeling for the problem,
one can try to do some numerical computations. However, this is not so easy either, since the CM values will be
large (or small) algebraic numbers in $F$. So from the floating point evaluation one cannot get the exact value.
However, using the result of \cite{BY} one obtains a convincing
heuristic how the values should look like. With this extra
information, the problem becomes accessible. For simplicity let us
assume that $m$ is square-free and that $p=5$ as before. (For $p=13,17$
the same argument applies.)


The value of the Petersson
metric of $H_m$ at $\calC\calM(K)$ can be computed by means of the
formula of \cite{BY}.  Let $\frakl\subset\calO_F$ be a prime ideal
above a prime $l\in \Z$. Using the notation of \cite{BY}, to $m$ and
$l$ we can associate the quantity
\begin{align}
h_m(l)=\frac{W_{\tilde K}}{4}(b_{mp^2}(l)-b_{m}(l)).
\end{align}
As explained in \cite{BY} (1.10) it should have a geometric interpretation as the intersection number of suitable
models of $\dv(H_m)$ and $\calC\calM(K)$ in the fiber above $l$ of the moduli space of abelian surfaces with
$\calO_F$-action and $\partial_F^{-1}$-polarization. In the same way, over $\calO_F$, we could associate to $m$ and
$\frakl$ the intersection number $h_m(\frakl)$ in the fiber above $\frakl$. We should have
\[
h_m(l)=\begin{cases}  2 h_m(\frakl),& \text{if $l$ is ramified in $\calO_F$,}\\
h_m(\frakl),& \text{if $l$ is inert in $\calO_F$,}\\
h_m(\frakl)+h_m(\frakl'),& \text{if $l$ is split in $\calO_F$.}
\end{cases}
\]
According to \eqref{eq:1}, one should be able to write $h_m(\frakl)=h_m^+(\frakl)+h_m^-(\frakl) $, where
$h_m^\pm(\frakl)$ denotes the intersection of $\dv(\Psi^\pm_m)$ and $\calC\calM(K)$ in the fiber above $\frakl$.
Since $(\Psi^+_m)'=\Psi^-_m$, we would have
\[
h_m^+(\frakl)=h_m^-(\frakl').
\]
In view of \eqref{eq:2} the intersection number of $\dv(\tilde\Psi_m)$ and $\calC\calM(K)$ in the fiber above
$\frakl$, would be given by
\begin{equation} \label{eq:3}
h_m^+(\frakl)-h_m^-(\frakl)=h_m^+(\frakl)-h_m^+(\frakl').
\end{equation}
But since $\tilde \Psi_m$ has weight $0$, this quantity would be equal to
$\ord_\frakl(\tilde\Psi_m(\calC\calM(K)))$. Consequently,
\begin{align}
 |\ord_\frakl(\tilde\Psi_m(\calC\calM(K)))|\leq
\begin{cases}
0,& \text{if $l$ is ramified or inert in $\calO_F$,}\\
h_m(l),& \text{if $l$ is split in $\calO_F$.}
\end{cases}
\end{align}
The quantities $h_m(l)$ can be computed by means of the formula of
\cite{BY}. In that way there are only finitely many possibilities left
for the prime ideal factorization of
$\tilde\Psi_m(\calC\calM(K))$. Using a computer algebra system one can
now compute the CM value if the Fourier expansion of $\tilde \Psi_m$
is known.

In the special cases $p=5$ and $m=1,2$ considered before, we computed a few CM values corresponding to non-biquadratic CM
fields $K$. We listed some data on a few CM fields in Table \ref{cmfielddata}, including the class number $h_K$, and a system of representatives for the ideal class group of $K$.
The corresponding CM values are given in Tables \ref{twistedvalues} and \ref{twistedvalues2}.

\begin{table}[h]
\caption{\label{twistedvalues} CM values of $\tilde\Psi_1$ for $\Q(\sqrt{5})$}
\begin{tabular*}{16.2cm}{|r|@{\extracolsep{\fill}}p{7cm}|l| }
\hline
\rule[-3mm]{0mm}{8mm}$q$ & $ \prod_l l^{h_1(l)}$ & $\tilde\Psi_1(\CM(K))$  \\
\hline
\rule[-3mm]{0mm}{8mm} $5$ & \scriptsize  $5^{40}\cdot 19^{10}\cdot 29^{10}$ &
$\left(\frac{\omega}{\omega'}\right)^{54}
\frac{(4+\omega')^{10}\cdot(5+\omega')^{10}}{(4+\omega)^{10}\cdot(5+\omega)^{10}}$ \\
\hline
\rule[-3mm]{0mm}{8mm}  $41$  &\scriptsize
$2^{16}\cdot 5^{42}\cdot 23^4\cdot 31^4\cdot 37^2\cdot 61^4\cdot 107^2\cdot 127^2$&

$\left(\frac{\omega}{\omega'}\right)^{32}\frac{(5+2\omega)^2\cdot (7+3\omega')^2}{(5+2\omega')^2\cdot (7+3\omega)^2}$\\
\hline
\rule[-3mm]{0mm}{8mm} $61$ & \scriptsize
$3^{12}\!\cdot\! 5^{42}\!\cdot\! 13^4\!\cdot\! 41^4\!\cdot \!83^2\! \cdot\! 103^2\!\cdot\! 109^2\!\cdot\! 113^2\!\cdot \!199^2\!\cdot\! 379^2 $&
$\left(\frac{\omega}{\omega'}\right)^{24}
\frac{(6+\omega)^2\cdot
                       (10+\omega')^2\cdot (13+3\omega)^2\cdot(19+\omega')^2}{
                       (6+\omega')^2\cdot (10+\omega)^2\cdot (13+3\omega')^2
                       \cdot (19+\omega)^2}$\\
\hline
\rule[-3mm]{0mm}{8mm} $109$ &\scriptsize $3^{12}\!\cdot\! 5^{40}\!\cdot \!7^6\!\cdot\! 43^2\!\cdot\! 61^4\!\cdot\! 71^4\!\cdot\! 73^2\!\cdot\! 97^2\!\cdot\! 113^2\!\cdot\! 223^2\!\cdot\! 409^2\!\cdot\! 499^2$&
$\left(\frac{\omega}{\omega'}\right)^{12}
\frac{(7+3\omega)^2\cdot (8+\omega)^2\cdot (19+3\omega)^2\cdot (20+9\omega)^2}{
(7+3\omega')^2\cdot(8+\omega')^2\cdot (19+3\omega')^2\cdot

(20+9\omega')^2}$
  \\
\hline
\rule[-3mm]{0mm}{8mm} $241$ &\scriptsize
$2^{50}\!\cdot\! 3^{36}\!\cdot\! 5^{126}\!\cdot\! 29^{14}\!\cdot\! 47^6\!\cdot\! 53^6\!\cdot\! 61^8\!\cdot\! 67^2\!\cdot\! 83^2\!\cdot\! 97^6\!\cdot\! 229^4\!\cdot\! 257^2\!\cdot\! 331^4\!\cdot\! 347^2\!\cdot\! 617^2$&
$\left(\frac{\omega}{\omega'}\right)^{84}
\frac{(5+\omega')^2\cdot (17+3\omega')^2}{(5+\omega)^2\cdot (17+3\omega)^2}$
  \\
\hline
\rule[-3mm]{0mm}{8mm} $281$ &\scriptsize
$2^{44}\!\cdot\! 5^{128}\!\cdot\! 7^{22}\!\cdot\! 17^{12}\!\cdot\! 43^6\!\cdot\! 53^6\!\cdot\! 59^6\!\cdot\! 101^4\!\cdot\! 109^6\!\cdot\! 137^2
\!\cdot\! 191^4\!\cdot\! 317^2\!\cdot\! 421^4\!\cdot\! 647^2\!\cdot\! 787^2\!\cdot\! 857^2\!\cdot\! 877^2$&
$\left(\frac{\omega}{\omega'}\right)^{72}
\frac{(7+2\omega)^2\cdot (10+\omega)^2\cdot (13+2\omega)^2\cdot(19+4\omega)^2}
{(7+2\omega')^2\cdot (10+\omega')^2\cdot (13+2\omega')^2\cdot (19+4\omega')^2}$
  \\
\hline
\rule[-3mm]{0mm}{8mm} $409$ &\scriptsize
$2^{48}\!\cdot\! 3^{42}\!\cdot\! 5^{120}\!\cdot\! 17^{12}\!\cdot\! 23^{14}\!\cdot\!
53^{10}\!\cdot\! 83^4\!\cdot\! 103^6\!\cdot\! 109^6\!\cdot\! 167^2\!\cdot\! 179^4\!\cdot\! 197^2\!\cdot\! 239^4\!\cdot\! 349^4
\!\cdot\! 571^4\!\cdot\! 1187^2\!\cdot\! 1277^2$&
$\left(\frac{\omega}{\omega'}\right)^{72}
\frac{(10+\omega)^2\cdot (23+2\omega)^2}
{(10+\omega')^{2}\cdot (23+2\omega')^2}$
  \\
\hline
\end{tabular*}
\end{table}

\begin{table}[h]
\caption{\label{twistedvalues2} CM values of $\tilde\Psi_2$ for $\Q(\sqrt{5})$}
\begin{tabular*}{16.2cm}{|r|@{\extracolsep{\fill}}p{5cm}|l| }
\hline
\rule[-3mm]{0mm}{8mm}$q$ & $ \prod_l l^{h_2(l)}$ & $\tilde\Psi_2(\CM(K))$  \\
\hline
\rule[-3mm]{0mm}{8mm} $5$ & \scriptsize
$2^{120}\cdot 5^{50}\cdot 89^{10}\cdot 109^{10}$ &
$\left(\frac{\omega}{\omega'}\right)^{102}
\frac{(9+\omega')^{10}\cdot(10+\omega')^{10}}{(9+\omega)^{10}\cdot(10+\omega)^{10}}
$\\
\hline
\rule[-3mm]{0mm}{8mm}  $41$  &\scriptsize
$2^{104}\!\cdot\! 5^{48}\!\cdot\! 23^{4}\!\cdot\! 37^2\!\cdot\! 41^2\!\cdot\! 43^2
\!\cdot \!73^2\!\cdot\!
83^2\!\cdot\! 113^2\!\cdot\! 349^2\!\cdot\! 449^2\!\cdot\! 769^2\!\cdot\! 829^2
\!\cdot \!1009^2$
&
$\left(\frac{\omega}{\omega'}\right)^{84}
\frac{ ( 6+  \omega)^2\cdot (17+5\omega)^2\cdot (19+8\omega)^2\cdot
(25+9\omega')^2\cdot (26+9\omega')^2\cdot (29+8\omega)^2}
{( 6+  \omega')^2\cdot (17+5\omega')^2\cdot (19+8\omega')^2\cdot
(25+9\omega)^2\cdot (26+9\omega)^2\cdot (29+8\omega')^2}$ \\
\hline
\rule[-3mm]{0mm}{8mm} $61$ & \scriptsize
$2^{120}\! \cdot \! 3^{14}\! \cdot\!  5^{48}\!  \cdot\!  13^2 \! \cdot\!  41^4\! \cdot\!  47^4\! \cdot\!  61^2\! \cdot \!
73^2\! \cdot \! 229^2\! \cdot\!  283^2\! \cdot\!  443^2\! \cdot\!  503^2\! \cdot\!  1489^2$
&
$\left(\frac{\omega}{\omega'}\right)^{80}
\frac{( 6+  \omega')^2\cdot (7+3\omega)^2\cdot (14+3\omega')^2\cdot
(35+11\omega)^2}{( 6+  \omega)^2\cdot (7+3\omega')^2\cdot (14+3\omega)^2
\cdot (35+11\omega')^2}
$\\
\hline
\rule[-3mm]{0mm}{8mm} $109$ &\scriptsize
$
2^{120}\!\cdot\! 3^{12}\! \cdot\!  5^{50}\!  \cdot\!  7^{12}\! \cdot\!  43^2\!
\cdot\!  263^2\! \cdot\!  281^4
\! \cdot\!  307^2\! \cdot\!  523^2\!  \cdot\!  683^2\! \cdot\!  823^2\!
\cdot\!  1429^2\! \cdot\!  2689^2$
&
$\left(\frac{\omega}{\omega'}\right)^{66}
 \frac{( 15+ 7\omega)^2\cdot (34+13\omega')^2\cdot (47+15\omega)^2}
{( 15+ 7\omega')^2\cdot (34+13\omega)^2\cdot (47+15\omega')^2}
$
\\
\hline
\end{tabular*}
\end{table}

\begin{table}[h]
\caption{\label{cmfielddata} CM extensions of $\Q(\sqrt{5})$}
\begin{tabular}{|r|l|l|l|}
\hline
\rule[-3mm]{0mm}{8mm}
$q$ & $K=F(\sqrt{\Delta})$ & $h_K$ & $\CL(K)$  \\
\hline
\rule[-3mm]{0mm}{8mm}
$5$ &$\Delta=-\frac{5+\sqrt{5}}{2}$& 1& $\mathcal{O}_K = \darstell{}{\sqrt{\Delta}}$\\
\hline
\rule[-3mm]{0mm}{8mm}
41 & $\Delta=- \frac{13+\sqrt{5}}{2}$ & 1& $\mathcal{O}_K =\darstell{}{\frac{1}{2}(\sqrt{\Delta}+\frac{3+\sqrt{5}}{2})}$\\
\hline
\rule[-3mm]{0mm}{8mm}
61 & $\Delta = - (9+2 \sqrt{5})$ & 1& $\mathcal{O}_K =\darstell{}{\frac{1}{2}(\sqrt{\Delta}+1)}$\\
\hline
\rule[-3mm]{0mm}{8mm}
109 & $\Delta = -\frac{21+\sqrt{5}}{2} $& 1& $\mathcal{O}_K =\darstell{}{\frac{1}{2}(\sqrt{\Delta}+\frac{3+\sqrt{5}}{2})}$\\
\hline
\rule[-3mm]{0mm}{8mm}
241 & $\Delta=- \frac{33+5 \sqrt{5}}{2}$& 3& $\mathcal{O}_K = \darstell{}{\frac{1}{2}(\sqrt{\Delta}+\frac{3+\sqrt{5}}{2})}$, \\
\rule[-3mm]{0mm}{8mm}&&&
$\mathfrak{A} =\darstell{2}{\frac{1}{2}(\sqrt{\Delta}+\frac{9+3 \sqrt{5}}{2})}$,\\
\rule[-3mm]{0mm}{8mm}&&&
$\mathfrak{B} =\darstell{4}{\frac{1}{2}(\sqrt{\Delta}+\frac{9+3 \sqrt{5}}{2})}$\\
\hline
\rule[-3mm]{0mm}{8mm}
281 & $\Delta=- \frac{37+7 \sqrt{5}}{2}$ & 3 & $\mathcal{O}_K = \darstell{}{\frac{1}{2}(\sqrt{\Delta}+\frac{1+\sqrt{5}}{2})}$,\\
\rule[-3mm]{0mm}{8mm}&&&
$\mathfrak{A} = \darstell{2}{\frac{1}{2}(\sqrt{\Delta}+\frac{1+\sqrt{5}}{2})}$,\\
\rule[-3mm]{0mm}{8mm}&&&
$\mathfrak{B} = \darstell{4}{\frac{1}{2}(\sqrt{\Delta}+\frac{9+\sqrt{5}}{2})}$\\
\hline
\rule[-3mm]{0mm}{8mm}
409 & $\Delta=- \frac{41+3 \sqrt{5}}{2}$ & 3&
$\mathcal{O}_K = \darstell{}{\frac{1}{2}(\sqrt{\Delta}+\frac{1+\sqrt{5}}{2})}$,\\
\rule[-3mm]{0mm}{8mm}&&&
$\mathfrak{A} = \darstell{2}{\frac{1}{2}(\sqrt{\Delta}+\frac{7+ 3 \sqrt{5}}{2})}$,\\
\rule[-3mm]{0mm}{8mm}&&&
$\mathfrak{B} = \darstell{4}{\frac{1}{2}(\sqrt{\Delta}+\frac{-1+3 \sqrt{5}}{2})}$\\
\hline
\end{tabular}
\end{table}

Combining the above considerations with Corollary 1.3 of \cite{BY} we are lead to the following conjecture.

\begin{conjecture}\label{conj}
Let $m$ be a positive integer, and let $\frakl\subset\calO_F$ be
a prime ideal above a prime $l\in \Z$ such that
$\ord_\frakl(\tilde\Psi_m(\calC\calM(K)))\neq 0$.
Then $\frakl'\neq \frakl$ and $4 l| m^2 p^2 q - r^2$ for some $r\in \Z$ with $|r|<mp\sqrt{q}$.
\end{conjecture}


\section{Further remarks and open problems}
\label{sect:7}

1. Is there a nice moduli interpretation of $\tilde \Psi_1$ as an
   invariant of abelian surfaces with $\calO_F$-multiplication and
   $\partial_F^{-1}$-polarization?

2. In the numerical calculations it always happened that
\[\tilde\Psi_m(\CM(K,\Phi,\calO_F))= \tilde\Psi_m(\CM(K,\sigma^3\Phi,\calO_F)).\]
This explains that the CM values in Tables \ref{twistedvalues} and
\ref{twistedvalues2} are squares. Is it possible to prove this in
general? (See also Remark~\ref{rem:square}.) Moreover, it is
striking that the CM values $\tilde\Psi_1(\CM(K,\Phi,\calO_F))$
corresponding to non-Galois CM fields are actually square-free
elements of $F$. Is this a coincidence or a general phenomenon?
Notice that by Corollary~\ref{clcor1.3} this would imply that
$\tilde\Psi_1(z)$ generates the class field $L_K$. This is the
case with the examples computed above.

3. Is there a (finite) subset $B$ of the positive integers such that the functions $\tilde \Psi_m$ ($m\in B$)
generate the function field of the symmetric Hilbert modular surface corresponding to $\Sl_2(\calO_F)$?

4. In terms of the generators $g_2$, $s_6$, $s_{10}$ of the ring
of symmetric Hilbert modular forms of even weight for
$\Q(\sqrt{5})$ the functions $\tilde \Psi_1$ and $\tilde \Psi_2$
look rather complicated. Are there other generators which yield a
nicer description (possibly of all $\tilde \Psi_m$)?

5. Describe how the correspondence of Theorem \ref{hilbert} behaves under the action of the corresponding Hecke
algebras. It should be Hecke equivariant, where the Hecke action on the image is multiplicative (see also \cite{Bo2} Problem 16.5).

6. Generalize the results of the present paper to Hilbert modular
surfaces of arbitrary discriminant.

\end{document}